\documentclass[11pt]{article}
\usepackage{amsmath}
\usepackage{multirow}
\usepackage{epsfig}

\def\inter{\mathop{\hbox{\rm int}}}

\def\Risk{{\mathop{\hbox{\rm Risk}}}}
 \def\ones{{\mathbf{1}}}
\def\W{{\cal W}}
\def\Ker{\mathop{\hbox{\rm Ker}}}
\def\V{{\cal V}}

 \def\sCard{\hbox{\rm\scriptsize Card}}
\newcommand{\wh}[1]{{\widehat{#1}}}
\usepackage[bf,normalsize,tableposition=top]{caption}
\usepackage{subfig}
\usepackage{amsfonts}
\usepackage{color}
\usepackage{epsfig}
\usepackage{srcltx}
\usepackage[colorlinks=true]{hyperref}
\def\barHplus{\hbox{\small$\left[\begin{array}{c|c}\bar{H}&\bar{h}\cr\hline\bar{h}^T&\cr\end{array}\right]$}}
\def\Z{{\cal Z}}
\def\F{{\cal F}}

\usepackage{amssymb}

\def\Diag{{\hbox{\rm Diag}}}
\usepackage{graphicx}
 \def\tid{\hbox{\rm\tiny id}}
\def\bR{{\mathbf{R}}}
\def\bZ{{\mathbf{Z}}}
\def\Opt{{\hbox{\rm Opt}}}

\def\A{{\cal A}}

\def\SG{{\cal SG}}

\def\cJ{{\cal J}}
\def\E{{\cal E}}
\def\P{{\cal P}}

\def\R{{\cal R}}
\def\G{{\cal G}}
\def\H{{\cal H}}
\def\M{{\cal M}}
\def\X{{\cal X}}

\def\N{{\cal N}}

\def\S{{\cal S}}

\def\bE{{\mathbf{E}}}
\def\bS{{\mathbf{S}}}
\def\Prob{\hbox{\rm Prob}}

\newtheorem{lemma}{Lemma}[section]
\newtheorem{proposition}{Proposition}[section]
\newtheorem{remark}{Remark}[section]
\def\Det{{\hbox{\rm Det}}}
\def\qed{\hfill$\Box$}
\def\Card{\mathop{\hbox{\rm Card}}}

\def\Tr{{\hbox{\rm Tr}}}

\newcommand{\be}{\begin{eqnarray}}
\newcommand{\ee}[1]{\label{#1}\end{eqnarray}}
\newcommand{\nn}{\nonumber \\}
\newcommand{\I}{{\cal I}}

\newcommand{\rf}[1]{~(\ref{#1})}
\newcommand{\half}{ \mbox{\small$\frac{1}{2}$}}

\newcommand{\ese}{\end{eqnarray*}}
\newcommand{\bse}{\begin{eqnarray*}}

\definecolor{MyDarkBlue}{rgb}{0,0.08,0.45}
\definecolor{MyViolet}{rgb}{0.45,0.08,0.95}
\definecolor{MyBrown}{rgb}{0.45,0.08,0}

\newcommand{\ai}[2]{{\color{blue}\  #2}}
\newcommand{\hide}[1]{}
\title{Estimating Linear and Quadratic forms via Indirect Observations}
\author{
Anatoli Juditsky
\thanks{LJK, Universit\'e Grenoble Alpes, 700 Avenue Centrale 38041 Domaine Universitaire
de Saint-Martin-d'H\`{e}res, France,
{\tt anatoli.juditsky@univ-grenoble-alpes.fr}}
\and Arkadi Nemirovski
\thanks{Georgia Institute
 of Technology, Atlanta, Georgia
30332, USA, {\tt nemirovs@isye.gatech.edu}\newline
The first author was supported by the LabEx PERSYVAL-Lab (ANR-11-LABX-0025) and the PGMO grant 2016-2032H. Research of
the second author was supported by NSF grants   CCF-1523768 and CMMI-1262063.}}
\date{}
\begin{document}
\maketitle
\begin{abstract}
In this paper, we further develop the approach, originating in {\cite{JN2009}}, to ``computation-friendly''  statistical estimation
via Convex Programming.{Our focus is} on estimating a linear or quadratic form of an unknown  ``signal,'' known to belong to a given convex compact set, via noisy indirect observations of the signal. {Classical} theoretical results on  the subject  deal with precisely stated statistical models and aim at designing statistical inferences and quantifying their performance in a closed analytic form.  In contrast to this traditional (highly instructive) descriptive framework, the approach
we promote here can be qualified as operational -- the estimation routines and their risks are {not available ``in a closed form,'' but are} yielded by an efficient computation. All we know in advance is that under favorable circumstances  the risk of the resulting estimate, whether high or low, is provably near-optimal under the circumstances. As a compensation for the lack of ``explanatory power,'' this approach
is applicable to a much wider family of observation schemes than those where ``closed form descriptive analysis'' is possible.
\par
{We discuss applications of this approach to classical problems of estimating linear forms of parameters of sub-Gaussian distribution and quadratic forms of partameters of Gaussian and discrete distributions. The performance of the constructed estimates is illustrated by computation experiments in which we compare the risks of the constructed estimates with (numerical) lower bounds for corresponding minimax risks for randomly sampled estimation problems.}
\end{abstract}

\section{Introduction}
 This paper can be considered as a follow-up to the paper \cite{PartI} dealing with hypothesis testing for {simple} families {-- families } of distributions specified in terms of upper bounds on their moment-generating functions. In what follows, we work with simple families of distributions, but our focus is on estimation of linear or quadratic forms of the unknown ``signal'' {(partly)} parameterizing the distribution in question. To give an impression of our approach and results, let us consider the sub-Gaussian case, where one is  given a random observation $\omega$ drawn from a sub-Gaussian distribution $P$ on $\bR^d$:
 $$
 \bE_{\omega\sim P}\{{\rm e}^{h^T\omega}\}\leq \mu^Th+\half h^T\Theta h\;\;\forall h\in\bR^d,
 $$
 with sub-Gaussianity parameters $\mu\in\bR^d$, $\Theta=\Theta^T\in\bR^{d\times d}$ affinely parameterized by ``signal'' $x\in\bR^m$. The goal is, given observation $\omega$ ``stemming'' from {\sl unknown} signal $x$ {\sl known to belong to a given convex compact set $X\subset\bR^m$}, to recover the value at $x$  of a given linear  form $g(\cdot):\bR^m\to\bR$.  The estimate $\widehat{g}$ we build
   is affine function of observation; the coefficients of the function, same as an upper bound on the $\epsilon$-risk of the estimate on $X$ \footnote{For the time being, given $\epsilon\in(0,1)$, $\epsilon$-risk of an estimate on $X$ is defined as the worst-case, over $x\in X$, width of $(1-\epsilon)$-confidence interval yielded by the estimate.} stem from  an optimal solution to an explicit convex optimization problem and thus can be specified in a computationally efficient fashion. Moreover, under  mild structural assumptions on the affine mapping $x\mapsto (\mu,\Theta)$ the resulting estimate is provably near-optimal in the minimax sense (see Section \ref{estlinformA} for details). The latter statement is an extension of the fundamental result of D. Donoho \cite{Don95} on near-optimality of affine recovery of a linear form  of signal in Gaussian observation scheme.
 \par
 This paper contributes to a long line of research on estimating linear (see, e.g., \cite{levit1975conditional,ibragimov1988estimation,efromovich1994adaptive,lepski1997optimal,klemela2001sharp,JN2009,butucea2009adaptive} and references therein) and quadratic  (\cite{hasminskii1980some,ibragimov1987some,bickel1988estimating,fan1991estimation,donu1990,
birge1995estimation,efromovich1996optimal,
laurent1997estimation,gayraud1999wavelet, huang1999nonparametric,laurent2000adaptive,
laurent2005adaptive,klemela2006sharp,
butmez2011} among others) functionals of parameters of probability distributions via observations drawn from these distributions.  In the majority of cited papers, the objective is to provide ``closed analytical form'' lower risk bounds for problems at hand and upper risk bounds for the proposed estimates, in good cases matching the lower bounds. This paradigm can be  referred to as ``descriptive;'' it relies upon analytical risk analysis and estimate design and possesses strong
explanation power. It, however, imposes severe restrictions on the structure of the statistical model, restrictions making the estimation problem amenable to  complete analytical treatment. There exists another, ``operational,''  line of research, initiated by D. Donoho in \cite{Don95}. The spirit of the operational approach is perfectly well illustrated by the main result of \cite{Don95} stating that when recovering the linear form of unknown signal $x$ known to belong to a given convex compact set $X$ via indirect Gaussian  observation $\omega=Ax+\xi$, $\xi\sim \N(0,I)$, the worst-case, over $x\in X$, risk of an affine in $\omega$ estimate yielded by optimal solution to an explicit convex optimization problem is  within the factor 1.2 of the minimax optimal risk. Subsequent ``operational'' literature is of similar spirit: both the recommended estimate and its risk are given by an efficient computation (typically, stem from solutions to explicit convex optimization problems); in addition, in good situations we know in advance  that the resulting risk, whether large or small, is nearly minimax optimal. The explanation power of  operational results is almost nonexisting; as a compensation, the scope of operational results  is usually much wider {than the one of analytical results}. For example, the just cited result of D. Donoho imposes no restrictions on $A$ and $X$, except for convexity and compactness of $X$; in contrast, all known {to us} analytical results on the same problem subject $(A,X)$ to severe structural restrictions. In terms of the outlined ``descriptive -- operational'' dichotomy, our paper is operational. For instance, in the problem of estimating
linear functional of signal $x$ affinely parameterising the parameters $\mu,\Theta$ of sub-Gaussian distribution    we started with, we allow for quite general affine mapping $x\to (\mu,\Theta)$ and for general enough signal set $X$, the only  restrictions on $X$ being convexity and compactness.\par
Technically, the approach we use in this paper combines the machinery developed in \cite{GJN2015,PartI} and the Cramer-type techniques for upper-bounding the risk of an affine estimate developed in \cite{JN2009}.\footnote{To handle the case of estimates quadratic in observation, we treat them as affine functions of ``quadratic lifting'' $\omega^+=[\omega;1][\omega;1]^T$ of the actual observation $\omega$.} On the other hand, this approach can also be viewed as ``computation-friendly'' extension of theoretical results on ``Cramer tests'' supplied by \cite{Birge1982,Birge1981,Birge1983,Birge2006} in conjunction with techniques of
 \cite{donoho1990minimax,donu1990,donoho1991geom2,donoho1991geom3,Don95,butmez2011}, which exploits the  most attractive, in our opinion, feature of this {line of} research
 -- potential applicability to a wide variety of observation schemes and (convex) signal sets $X$.
 \par
 The rest of the paper is organized as follows. In Section \ref{sect:Setup} we, following \cite{PartI}, describe the families of distributions we are working with. We present the estimate construction and study its general properties in Section \ref{estlinform}. Then in Section \ref{estlinformA} we discuss {applications} to {estimating} linear forms of sub-Gaussian distributions. In Section \ref{sec:quadratic} we apply the proposed construction  to estimating quadratic forms of parameters of Gaussian and discrete distributions. To illustrate the performance of the proposed approach we describe results of some preliminary numerical experiments in which we compare the bounds
 on the risk of estimates supplied by our machinery with (numerically computed) lower bounds on the minimax risk.
To streamline the presentation, all proofs are collected in the appendix.
\bigskip\par\noindent
{\bf Notation.} In what follows, $\bR^n$ and $\bS^n$ stand for the spaces of real $n$-dimensional vectors and real symmetric $n\times n$ matrices, respectively; both spaces are equipped with the standard inner products, $x^Ty$, resp., $\Tr(XY)$.  Relation $A\succeq B$ ($A\succ B$) means that $A$, $B$ are symmetric matrices of the same size such that $A-B$ is positive semidefinite (resp., positive definite). We denote $\bS^n_+=\{S\in \bS^n:S\succeq0\}$ and $\bS^n_{++}=\inter\bS^n_+=\{S\in\bS^n:S\succ0\}$.\par
We use ``MATLAB notation:'' $[X_1;...;X_k]$ means vertical concatenation of matrices $X_1,...,X_k$ of the same width, and $[X_1,...,X_k]$ means  horizontal concatenation of matrices $X_1,...,X_k$ of the same height. In particular, for reals $x_1,...,x_k$, $[x_1;...;x_k]$ is a $k$-dimensional column vector with entries $x_1,...,x_k$.\par
For probability distributions $P_1,...,P_K$, $P_1\times...\times P_K$ is the product distribution on the direct product of the corresponding probability spaces; when $P_1=...=P_K$, we denote $P_1\times...\times P_K$ by $P^K$ or $[P]^K$.
\par
Given positive integer $d$, $\theta\in \bR^d$, $\Theta
\in
\bS^d_+$, we denote by $\SG(\theta,\Theta)$ the family of all sub-Gaussian, with parameters $(\theta,\Theta)$, probability distributions, that is, the family of all Borel probability distributions $P$ on $\bR^d$ such that
$$
\forall f\in\bR^d: \ln\left(\bE_{\zeta\sim P}\{\exp\{f^T\zeta\}\}\right)\leq f^T\theta+{1\over 2}f^T\Theta f.
$$
We use shorthand notation $\omega\sim\SG(\theta,\Theta)$ to express the fact that the probability distribution of random vector $\omega$ belongs to the family $\SG(\theta,\Theta)$.

\section{Simple families of probability distributions}\label{sect:Setup}
Let
\begin{itemize}
\item $\F$, $0\in\inter \F$, be a closed convex set in $\Omega=\bR^m$ symmetric w.r.t. the origin,
\item $\M$ be a closed convex set in some $\bR^n$,
\item $\Phi(h;\mu):\F\times\M \to\bR$ be a continuous function convex in $h\in \F$ and concave in $\mu\in \M$.
\end{itemize}
Following \cite{PartI}, we refer to $\F,\M,\Phi(\cdot,\cdot)$ satisfying the above restrictions as to {\sl regular data}. Regular data
$\F,\M,\Phi(\cdot,\cdot)$ define the family
$$
\S=\S[\F,\M,\Phi]
$$
of Borel probability distributions $P$ on $\Omega$ such that

\begin{equation}\label{eq1a}
\begin{array}{l}
\exists \mu\in\M: \forall h\in\F: \ln\left(\int_\Omega\exp\{h^T\omega\}P(d\omega)\right)\leq \Phi(h;\mu).
\end{array}
\end{equation}
We say that distributions satisfying \rf{eq1a} are {\em simple}. Given regular data $\F,\M,\Phi(\cdot,\cdot)$, we refer to $\S[\F,\M,\Phi]$ as to {\sl simple} family of distributions associated with the data $\F$, $\M$, $\Phi$.
{Standard examples of simple families are supplied by ``good observation schemes,'' as defined in \cite{JN2009,GJN2015}, and include the families of Gaussian, Poisson and discrete distributions. For other instructive} examples  and an algorithmic ``calculus'' of simple families, the reader is referred to \cite{PartI}. {We present here three examples of simple families which we use in the sequel.}
\subsection{Sub-Gaussian distributions}\label{sect:subGauss}
Let $\F=\Omega=\bR^d$, $\M$ be a closed convex subset of the set $\G_d=\{\mu=(\theta,\Theta):\theta\in\bR^d,\Theta\in\bS^d_+\}$,  and let
$$
\Phi(h;\theta,\Theta)=\theta^Th+\half h^T\Theta h.
$$
In this case, $\S[\F,\M,\Phi]$ contains all sub-Gaussian distributions $P$ on $\bR^m$ with sub-Gaussianity parameters from $\M$:
\begin{equation}\label{eq40}
(\theta,\Theta)\in \M\Rightarrow \SG(\theta,\Theta)\subset\S[\F,\M,\Phi].
\end{equation}
In particular, $\S[\F,\M,\Phi]$ contains all Gaussian distributions $\N(\theta,\Theta)$ with $(\theta,\Theta)\in\M$.

\subsection{Quadratically lifted Gaussian observations}\label{sect:GaussLift}
Let $\V$ be a nonempty convex compact subset of $\bS^d_+$. This set gives rise to the family $\P_\V$ of distributions of {\sl quadratic liftings} $[\zeta;1][\zeta;1]^T$ of random vectors $\zeta\sim\N(\theta,\Theta)$ with $\theta\in\bR^d$ and $\Theta\in\V$. Our goal now is to build regular data such that the associated simple family of distributions contains $\P_\V$. To this end we select
$\Theta_*\in \bS^d_{++}$ and $\delta\geq 0$  such that for all $\Theta\in\V$ one has
\be
 \Theta\preceq \Theta_*, \;\;\hbox{and}\;\; \|\Theta^{1/2}\Theta_*^{-1/2}-I\|\leq \delta,
\ee{56delta}
where $\|\cdot\|$ is the spectral norm; under these restrictions, the smaller are $\Theta_*$ and $\delta$, the better.
Observe that for all $\Theta\in\V$, we have $0\preceq \Theta_*^{-1/2}\Theta\Theta_*^{-1/2}\preceq I$. Hence
\[
\|\Theta^{1/2}\Theta_*^{-1/2}\|^2=\|\Theta_*^{-1/2}\Theta^{1/2}\|^2=\|\Theta_*^{-1/2}\Theta^{1/2}[\Theta_*^{-1/2}\Theta^{1/2}]^T\|
=\|\Theta_*^{-1/2}\Theta\Theta_*^{-1/2}\|
\leq 1,
\]
and we lose nothing when assuming from now on that $\delta\in[0,2]$.
The required regular data are given by the following
\begin{proposition}\label{propGausslift0} In the just described situation, let $\gamma\in(0,1)$,
\[
\Z^+=\{Z\in \bS^{d+1}:\;Z_{d+1,d+1}=1\},\;\;\H_\gamma=\{H\in\bS^d: -\gamma\Theta_*^{-1}\preceq H\preceq \gamma\Theta_*^{-1}\}
\]
and let $\F=\bR^d\times\H_\gamma$, $\M^+=\V\times \Z^+$. We set
\be
\Phi(h,H;\Theta,Z)&=&\Upsilon(H,\Theta)+\Gamma(h,H,Z),\\
\Upsilon(H,\Theta)&=&-\half \ln\Det(I-\Theta_*^{1/2}H\Theta_*^{1/2})+\half \Tr([\Theta-\Theta_*]H)\nn&&+{\delta(2+\delta)\over 2(1-\|\Theta_*^{1/2}H\Theta_*^{1/2}\|)}\|\Theta_*^{1/2}H\Theta_*^{1/2}\|_F^2\nn
\Gamma(h,H;Z)&=&\half \Tr\left(Z\left[\hbox{\small$\left[\begin{array}{c|c}H&h\cr\hline h^T&\end{array}\right]$}+[H,h]^T[\Theta_*^{-1}-H]^{-1}[H,h]\right]\right).\nonumber
\ee{PhisubG}
Then
\item[(i)]  $\F,\M^+,\Phi$ form a regular data, and for every $(\theta,\Theta)\in\bR^d\times\V$ it holds for all $(h,H)\in\F$:
\begin{equation}\label{itholds123}
\ln\left(\bE_{\zeta\sim\N(\theta,\Theta)}\left\{{\rm e}^{h^T\zeta+{1\over 2}\zeta^TH\zeta}\right\}\right)\leq \Phi\left(h,H;\Theta,[\theta;1][\theta;1]^T\right).
\end{equation}
\item[(ii)] Besides this, function $\Phi(h,H;\Theta,Z)$ is coercive in the convex argument: whenever $(\Theta,Z)\in\M^+$, $(h_i,H_i)\in\F$ and $\|(h_i,H_i)\|\to\infty$ as $i\to\infty$, we have $\Phi(h_i,H_i;\Theta,Z)\to\infty$.
\end{proposition}
For proof, see Appendix \ref{sec:ppropGausslift}.
\subsection{Quadratically lifted discrete observations}\label{sect:DiscrLift}
Consider a random variable $\zeta\in \bR^d$ taking values $e_i$, $i=1,...,d$, where $e_i$ are standard basic orths in $\bR^d$.\footnote{This is nothing more than a convenient way of thinking of a discrete random variable taking values in a $d$-element set.} We identify the probability distribution $P_\mu$ of such variable with a point $\mu=[\mu_1;...;\mu_d]$ from the $d$-dimensional probabilistic simplex
\def\bDelta{{\mathbf{\Delta}}}
$
\Delta^d=\{\nu\in\bR^d_+:\sum_{i=1}^d \nu_i=1\}
$
where $\mu_i=\Prob\{\zeta=e_i\}$. Let now $\zeta^{K}=(\zeta_1,...,\zeta_{K})$ with $\zeta_k$ drawn independently across $k$ from  $P_\mu$, and let
\be
\omega[\zeta^K]={2\over K(K-1)}\sum_{1\leq j<j\leq K} \omega_{ij}[\zeta^K],\;\;\;\;
\omega_{ij}[\zeta^{K}]=\half[\zeta_i\zeta_j^T+\zeta_j\zeta_i^T],\,1\leq i<j\leq K.
\ee{liftedo1}
We are about to point our regular data such that the associated simple family of  distributions contains the distributions of the ``quadratic lifts'' $\omega[\zeta^K]$ of random vectors $\zeta^K$.
\begin{proposition}\label{discretelift}
Let  $\F=\bS^d$,
\be
\bDelta^d=\left\{Z\in \bS^d: \,Z_{ij}\geq 0\,\forall i,j,\, \sum_{i,j}Z_{ij}=1\right\}.
\ee{mdeltad}
and let $\Z^d$ be a set of all positive semidefinite matrices from $\bDelta^d$.
Denote
\be
\Phi(H;Z)=\ln\left(\sum_{i,j=1}^m Z_{ij}\exp\{H_{ij}\}\right):\;\bS^d\times \bDelta^d\to \bR,
\ee{quaddphi}
so that $\Phi(\cdot;\cdot)$ is convex-concave on $\bS^d\times \bDelta^d$.
We set
\[
\Phi_M(H;Z)=M\Phi(H/M;Z),\;\;M\in \bZ_+.
\]
Then for $M=M(K)=\lfloor K/2\rfloor$,
\be
\ln\left(\bE_{\zeta^K\sim P_\mu^K}\left\{\exp\{\Tr(H\omega[\zeta^K])\}\right\}\right)\leq \Phi_M(H;\mu\mu^T).
\ee{excovering}
In other words, the simple family $\S[\F,\Z^d,\Phi_{\lfloor K/2\rfloor}]$ contains distributions of all random variables $\omega[\zeta^K]$ with $\zeta\sim P_\mu$, $\mu\in \Delta^d$.
\end{proposition}
For proof, see Appendix \ref{sec:pprop2}.

\section{Estimating linear forms}\label{estlinform}
\subsection{Situation and goal}\label{linformsit}
Consider the situation as follows: given are Euclidean spaces $\E_F,$ $\E_M$, $\E_X$ along with
\begin{itemize}
\item regular data $\F\subset \E_F,\M\subset \E_M,\Phi(\cdot;\cdot):\F\times\M\to\bR$,
\item a nonempty set  $X$ contained in a convex compact set $\X\subset \E_X$,
\item an affine mapping $x\mapsto \A(x):\E_X\to\E_M $ such that $\A(\X)\subset\M$,
\item a vector $g\in \E_X$ and a constant $c$ specifying the linear form $G(x)=\langle g,x\rangle+c:\E_X\to\bR$ \footnote{from now on, $\langle u,v\rangle$ denotes the inner product of vectors $u,v$ belonging to a Euclidean space; what is this space, it always will be clear from the context.},
\item a tolerance $\epsilon\in(0,1)$.
\end{itemize}
Let $\P$ be the  family of all Borel probability distributions on $\E_F$.
Given a random observation
\begin{equation}\label{eqoservation}
\omega\sim P(\cdot)
\end{equation}
where $P\in\P$ is {\sl associated with} unknown signal  $x$ known to belong to $X$, ``association'' meaning that
\begin{equation}\label{123cond21}
\forall f\in\F: \ln\left(\int_{\E_F} {\rm e}^{\langle f,\omega\rangle}P(d\omega)\right)\leq \Phi(f;\A(x)),
\end{equation}
 we want to recover the quantity $G(x)$.
\par
Given $\rho>0$, we call an estimate -- a Borel function $\widehat{g}(\cdot):\E_F\to\bR$ -- $(\rho,\epsilon)$-accurate, if for all pairs $x\in X$, $P\in\P$ satisfying (\ref{123cond21}) it holds
\[
\Prob_{\omega\sim P}\left\{|\widehat{g}(\omega)-G(x)|>\rho\right\}\leq\epsilon.
\]
If $\rho_*$ is the infimum of those $\rho$ for which estimate $\widehat{g}$ is $(\rho,\epsilon)$-accurate, then clearly $\widehat{g}$ is
$(\rho_*,\epsilon)$-accurate. We refer to $\rho_*$ as the {\sl $\epsilon$-risk}
of the estimate $\widehat{g}$ w.r.t. the data $G(\cdot)$, $X$,  and $(\A,\F,\M,\Phi)$:
\begin{equation}\label{riskdefinition}
\begin{array}{rcl}
\Risk_\epsilon(\widehat{g}(\cdot)|G,X,\A,\F,\M,\Phi)&=&\min\bigg\{\rho:\forall (x,P)\in X\times\P:\\
&&\qquad\quad
\begin{array}{l}
\Prob_{\omega\sim P}\{\omega:\,|\widehat{g}(\omega)-G(x)|>\rho\}\leq\epsilon\\
\ln\left(\int{\rm e}^{\langle f,\omega\rangle}P(d\omega)\right)\leq\Phi(f;\A(x))\;\forall f\in\F
\end{array}\bigg\}\\
\end{array}
\end{equation}
When $G,X,\A,\F,\M,\Phi$ are clear from the context, we shorten $\Risk_\epsilon(\widehat{g}(\cdot)|G,X,\A,\F,\M,\Phi)$ to $\Risk_\epsilon(\widehat{g}(\cdot))$.
\par
{In the setting of} this section, we are about to build, in a computationally efficient fashion, an affine estimate $\widehat{g}(\omega)=\langle f_*,\omega\rangle+\varkappa$ along with $\rho_*$ such that  the estimate is $(\rho_*,\epsilon)$-accurate.
\subsection{The construction}\label{theconstruc}
 Let us set
$$
\F^+=\{(f,\alpha):f\in\E_F,\alpha>0,f/\alpha\in\F\}
$$
so that $\F^+$ is a nonempty convex set in $\E_F\times\bR_+$, and let
\[
\begin{array}{rcl}
\Psi_+(f,\alpha)&=&\sup\limits_{x\in \X}\left[\alpha\Phi(f/\alpha,\A(x))-G(x)\right]:\F^+\to\bR,\\
\Psi_-(f,\beta)&=&\sup\limits_{x\in \X}\left[\beta\Phi(-f/\beta,\A(x))+G(x) \right]:\F^+\to\bR,\\
\end{array}
\]
so that $\Psi_{\pm}$ are convex real-valued functions on $\F^+$ (recall that $\Phi$ is convex-concave and continuous on $\F\times\M$, while $\A(\X)$ is a compact subset of $\M$). These functions give rise to convex functions $\wh{\Psi}_{\pm}:\,\E_F\to\bR$ given by
$$\begin{array}{rcl}
\widehat{\Psi}_+(f)&:=&\inf_{\alpha}\left\{\Psi_+({f},{\alpha})+\alpha\ln(2/\epsilon):\alpha>0,(f,\alpha)\in\F^+\right\},\\
\widehat{\Psi}_-(f)&:=&\inf_{\alpha}\left\{\Psi_-({f},{\alpha})+\alpha\ln(2/\epsilon):\alpha>0,(f,\alpha)\in \F^+\right\}\\
\end{array}
$$
and to convex optimization problem
\begin{equation}\label{whynottoconsider}
\Opt=\min_f\left\{\widehat{\Psi}(f):=\half\left[\widehat{\Psi}_+(f)+\widehat{\Psi}_-(f)\right]\right\},
\end{equation}
With our approach, a ``presumably good'' estimate of $G(x)$ and its risk are given by an optimal (or nearly so) solution to the latter problem. The corresponding result is as follows:
\begin{proposition}\label{corinf}
In the situation  of Section \ref{linformsit}, let $\Phi$ satisfy the relation
\begin{equation}\label{Phisatisfies}
\Phi(0;\mu)\geq0\,\,\forall \mu\in\M.
\end{equation} Then
\be
\widehat{\Psi}_+(f)&:=&\inf_{\alpha}\left\{\Psi_+({f},{\alpha})+\alpha\ln(2/\epsilon):\alpha>0,(f,\alpha)\in\F^+\right\}\nn
&=&\max_{x\in\X}\inf_{\alpha>0,(f,\alpha)\in\F^+}\left[\alpha\Phi(f/\alpha,\A(x))-G(x)+\alpha\ln(2/\epsilon)\right],
\ee{psisareequala}
\be
\widehat{\Psi}_-(f)&:=&\inf_{\alpha}\left\{\Psi_-({f},{\alpha})+\alpha\ln(2/\epsilon):\alpha>0,(f,\alpha)\in \F^+\right\}\nn
&=&\max_{x\in\X}\inf_{\alpha>0,(f,\alpha)\in\F^+}\left[\alpha\Phi(-f/\alpha,\A(x))+G(x)+\alpha\ln(2/\epsilon)\right],
\ee{psisareequalb}
and the functions $\widehat{\Psi}_{\pm}(\cdot)$ are convex real-valued. {Furthermore,} a feasible solution
 $\bar{f}$,  $\bar{\varkappa}$, $\bar{\rho}$ to the system of convex constraints
\begin{equation}\label{sssuchthatcor}
\begin{array}{rcl}
\widehat{\Psi}_+(f)\leq{\rho}-{\varkappa},\;\;
\widehat{\Psi}_-(f)\leq{\rho}+{\varkappa}
\end{array}
\end{equation}
in variables $f$, $\rho$, $\varkappa$ induces estimate
\begin{equation}\label{neweqlin}
\widehat{g}(\omega)=\langle \bar{f},\omega\rangle+\bar{\varkappa},
\end{equation}
 of $G(x)$, $x\in X$,  with $\epsilon$-risk at most $\bar{\rho}$:
\begin{equation}\label{RiskBoundSep}
\Risk_\epsilon(\widehat{g}(\cdot)|G,X,\A,\F,\M,\Phi)\leq\bar{\rho}.
\end{equation}
Relation {\rm (\ref{sssuchthatcor})} (and thus -- the risk bound {\rm (\ref{RiskBoundSep})}) clearly holds true  when $\bar{f}$ is a candidate solution
to  problem  {\rm (\ref{whynottoconsider})} and
$$
\bar{\rho}=\widehat{\Psi}(\bar{f}),\,\,\,\bar{\varkappa}=\half\left[\widehat{\Psi}_-(\bar{f}) - \widehat{\Psi}_+(\bar{f})\right].
$$
As a result, by properly selecting $\bar{f}$ we can make (an upper bound on) the $\epsilon$-risk of estimate {\rm (\ref{neweqlin})}  arbitrarily close to $\Opt$, and equal to $\Opt$ when optimization problem {\rm (\ref{whynottoconsider})} is solvable.
\end{proposition}
For proof, see Appendix \ref{proof:corinf}.
\subsection{Estimation from repeated observations}\label{sectrepeated}
Assume that in the situation described in Section \ref{linformsit} we have access to $K$ observations $\omega_1,...,\omega_K$ sampled, independently of each other, from a probability distribution $P$, and are allowed to build our estimate based on these $K$ observations rather than on a single observation.
We can immediately reduce this new situation to the previous one simply by redefining the data. Specifically, given  $\F\subset \E_F$, $\M\subset\E_M$, $\Phi(\cdot;\cdot):\F\times\M\to\bR$, $X\subset \X\subset \E_X$, $\A(\cdot)$, $G(x)=\langle g,x\rangle
+c$, see Section \ref{linformsit}, and a positive integer $K$,
let us
 replace $\F\subset \E_\F$ with $\F^K:=\underbrace{\F\times...\times \F}_{K}\subset \E_\F^K:=\underbrace{\E_\F\times...\times \E_\F}_{K}$, and replace  $\Phi(\cdot,\cdot):\F\times\M\to\bR$ with $\Phi^K(f^K=(f_1,...,f_K);\mu)=\sum_{i=1}^K\Phi(f_i;\mu):\F^K\times\M\to\bR$. It is immediately seen that the updated data satisfy all requirements imposed on the data in Section \ref{linformsit}. Furthermore, for all
$f^K=(f_1,...,f_K)\in\F^K$,
whenever a Borel probability distribution $\P$ on $\E_\F$ and $x\in X$ satisfy (\ref{123cond21}), the distribution $P^K$ of $K$-element i.i.d. sample $\omega^K=(\omega_1,...,\omega_K)$ drawn from $P$ and $x$ are linked by the relation
\begin{equation}\label{123cond21repeated}
 \ln\left(\int_{\E_F^K} {\rm e}^{\langle f^K,\omega^K\rangle}P^K(d\omega^K)\right)
=\sum_i\ln\left(\int_{\E_F} {\rm e}^{\langle f_i,\omega_i\rangle}P(d\omega_i)\right)\leq \Phi^K(f^K;\A(x)).
\end{equation}
Applying to our new data the construction from Section \ref{theconstruc}, we arrive at ``repeated observations'' version of 
Proposition \ref{corinf}. Note that the resulting convex constraints/objectives are symmetric w.r.t. permutations of the components $f_1,...,f_K$ of $f^K$, implying that we lose nothing when restricting ourselves with collections $f^K$ with equal to each other components; it is convenient to denote the common value of these components $f/K$. With these observations, Proposition 
\ref{corinf}
becomes the statements as follows (we use the assumptions and the notation from the previous section):

\begin{proposition}\label{corinfrepeated}
In the situation described in Section \ref{linformsit}, let $\Phi$ satisfy the relation {\rm (\ref{Phisatisfies})}, and let a positive integer $K$ be given.
Then functions $\widehat{\Psi}_\pm:\E_F\to\bR$,
\[
\begin{array}{rcl}
\widehat{\Psi}_+(f)&:=&\inf_{\alpha}\left\{\Psi_+({f},{\alpha})+K^{-1}\alpha\ln(2/\epsilon):\alpha>0,(f,\alpha)\in\F^+\right\}\\
&=&\max_{x\in\X}\inf_{\alpha>0,(f,\alpha)\in\F^+}\left[\alpha\Phi(f/\alpha,\A(x))-G(x)+K^{-1}\alpha\ln(2/\epsilon)\right],
\\
\widehat{\Psi}_-(f)&:=&\inf_{\alpha}\left\{\Psi_-({f},{\alpha})+K^{-1}\alpha\ln(2/\epsilon):\alpha>0,(f,\alpha)\in \F^+\right\}\\
&=&\max_{x\in\X}\inf_{\alpha>0,(f,\alpha)\in\F^+}\left[\alpha\Phi(-f/\alpha,\A(x))+G(x)+K^{-1}\alpha\ln(2/\epsilon)\right]
\end{array}
\]
are convex and real valued.
Furthermore, let $\bar{f}$,  $\bar{\varkappa}$, $\bar{\rho}$ be a feasible solution to the system of convex constraints
\begin{equation}\label{sssuchthatcorrepeated}
\begin{array}{rcl}
\widehat{\Psi}_+(f)\leq{\rho}-{\varkappa},\;\;\;
\widehat{\Psi}_-(f)\leq{\rho}+{\varkappa}\\
\end{array}
\end{equation}
in variables $f$, $\rho$, $\varkappa$.
Then, setting
\begin{equation}\label{neweqlin1}
\widehat{g}(\omega^K)=\left\langle \bar{f},{1\over K}{\sum}_{i=1}^K\omega_i\right\rangle+\bar{\varkappa},
\end{equation}
we get an estimate of $G(x)$, $x\in X$, via independent $K$-repeated observations
$$
\omega_i\sim P,\,i=1,...,K
$$ with $\epsilon$-risk at most $\bar{\rho}$, meaning that whenever a Borel probability distribution $P$ is associated with  $x\in X$ in the sense of  {\rm(\ref{123cond21})}, one has
\[
\Prob_{\omega^K\sim P^K}\left\{\omega^K:|\widehat{g}(\omega^K)-G(x)|>\bar{\rho}\right\}\leq\epsilon.
\]

Relation {\rm (\ref{sssuchthatcorrepeated})}  clearly holds true  when $\bar{f}$ is a candidate solution to the convex optimization problem
\begin{equation}\label{whynottoconsiderrepeated}
\Opt=\min_f\left\{\widehat{\Psi}(f):=\half\left[\widehat{\Psi}_+(f)+\widehat{\Psi}_-(f)\right]\right\}
\end{equation}
and
$$
\bar{\rho}=\widehat{\Psi}(\bar{f}),\,\,\bar{\varkappa}=\half \left[{\widehat{\Psi}_-(\bar{f}) - \widehat{\Psi}_+(\bar{f})}\right].
$$
As a result, properly selecting $\bar{f}$, we can make (an upper bound on) the $\epsilon$-risk of estimate $\widehat{g}(\cdot)$  arbitrarily close to $\Opt$, and equal to $\Opt$ when optimization problem {\rm (\ref{whynottoconsiderrepeated})} is solvable.
\end{proposition}
From now on, if otherwise is not explicitly stated, we deal with $K$-repeated observations; to get back to single-observation case, it suffices to set $K=1$.

\section{Application: estimating linear form of  parameters {of sub-Gaussian distributions}}\label{estlinformA}
\subsection{Situation}
We are about to apply construction form Section \ref{estlinform} in the situation where our observation is sub-Gaussian with parameters affinely parameterized by signal $x$, and our goal is to recover a linear function of $x$. Specifically, consider the   situation described in  Section  \ref{estlinform}, with the data as follows:
\begin{itemize}
\item $\F=\E_F=\bR^d$, $\M=\E_M=\bR^d\times\bS_+^d$, $\Phi(h;\mu,M)=h^T\mu+{1\over 2}h^TMh:\bR^d\times(\bR^d\times\bS_+^d)\to\bR$ (so that $\S[\F,\M,\Phi]$ is the family of all sub-Gaussian distributions on $\bR^d$);
\item $X=\X\subset\E_X=\bR^{n_x}$ is a nonempty convex compact set, and
\item $\A(x)=(Ax+a,M(x))$, where $A$ is $d\times n_x$  matrix, and $M(x)$ is affinely depending on $x$ symmetric $d\times d$ matrix such that $M(x)$ is $\succeq0$ when $x\in X$,
\item $G(x)$ is an affine function on $\E_X$.
\end{itemize}
Same as in Section \ref{estlinform}, our goal is to recover the value of a given linear function $G(y)=g^Ty+c$ at unknown signal $x\in \X$ via $K$-repeated observation $\omega^K=(\omega_1,...,\omega_K)$ with
$\omega_i$ drawn, independently across $i$, from a distribution $P$ which is associated with $x$, which now means ``is sub-Gaussian with parameters $(Ax+a,M(x))$.'' We refer to {\sl Gaussian case} as to the special case of the just described problem, where the distribution $P$ associated with signal $x$ is exactly $\N(Ax+a,M(x))$.
\par
In the case in question $\Phi(0;\mu,M)=0$, so that (\ref{Phisatisfies}) takes place, and the left hand sides in the constraints (\ref{sssuchthatcorrepeated}) are
$$
\begin{array}{rcl}
\widehat{\Psi}_+(f)&=&\sup_{x\in X}\inf_{\alpha>0}
\left\{f^T[Ax+a]+{1\over 2\alpha} f^TM(x)f+K^{-1}\alpha\ln(2/\epsilon)- G(x)\right\}\\
&=&\max_{x\in X} \left\{
{\left[2K^{-1}\ln(2/\epsilon)f^TM(x)f\right]^{1/2}}+f^T[Ax+a]-G(x)\right\},\\
\widehat{\Psi}_-(f)&=&\sup_{x\in X}\inf_{\alpha>0}
\left\{-f^T[Ax+a]+{1\over 2\alpha} f^TM(x)f+K^{-1}\alpha\ln(2/\epsilon)+ G(x)\right\}\\
&=&\max_{x\in X} \left\{
{\left[2K^{-1}\ln(2/\epsilon)f^TM(x)f\right]^{1/2}}-f^T[Ax+a]+G(x)\right\}.\\
\end{array}
$$
Thus,  system (\ref{sssuchthatcorrepeated}) reads
$$
\begin{array}{rcl}
a^Tf+\max\limits_{x\in X} \left\{
\left[2K^{-1}\ln(2/\epsilon)f^TM(x)f\right]^{1/2}+f^TAx-G(x)\right\}&\leq& \rho-\varkappa,\\
-a^Tf+\max\limits_{x\in X} \left\{
\left[2K^{-1}\ln(2/\epsilon)f^TM(x)f\right]^{1/2}-f^TAx+G(x)\right\}&\leq& \rho+\varkappa.\\
\end{array}
$$

We arrive at the following version of Proposition \ref{corinfrepeated}:
\begin{proposition}\label{propwearrive} In the situation  described {above,} given $\epsilon\in(0,1)$, let $\bar{f}$ be a feasible solution to the convex optimization problem
\begin{equation}\label{bestprob}
\Opt=\min_{f\in\bR^d}\left\{\widehat{\Psi}(f):=\half \left[\widehat{\Psi}_+(f)+\widehat{\Psi}_-(f)\right]\right\}
\end{equation}
where
\bse
\widehat{\Psi}_+(f)&=&
\max_{x\in X}\left\{
\left[2K^{-1}\ln(2/\epsilon)f^TM(x)f\right]^{1/2}+f^TAx-G(x)\right\}+a^Tf,\\
\widehat{\Psi}_-(f)&=&\max_{y\in X}\left\{
\left[2K^{-1}\ln(2/\epsilon)f^TM(x)f\right]^{1/2} -f^TAy+G(y)\right\}-a^Tf.
\ese
Let us set
\[
\bar{\varkappa}=\half\left[\widehat{\Psi}_-(\bar{f})-\widehat{\Psi}_+(\bar{f})\right],\;\;\bar{\rho}=\widehat{\Psi}(\bar{f}).
\]
Then the $\epsilon$-risk of the affine estimate
\[\widehat{g}(\omega^K)={1\over K}\sum_{i=1}^K\bar{f}^T\omega_i+\bar{\varkappa},
\]
taken w.r.t. the data listed in the beginning of this section, is at most $\bar{\rho}$.
\end{proposition}
It is immediately seen that optimization problem (\ref{bestprob}) is solvable, provided that $\bigcap\limits_{x\in X}\Ker(M(x))=\{0\}$, and an optimal solution $f_*$ to the problem, taken along with
\begin{equation}\label{takenalong}
\varkappa_*=\half\left[\widehat{\Psi}_-(f_*)-\widehat{\Psi}_+(f_*)\right],
\end{equation}
yields the affine estimate
$$
\widehat{g}_*(\omega)={1\over K}\sum_{i=1}^Kf_*^T\omega_i+\varkappa_*
$$
with $\epsilon$-risk, w.r.t. the data listed in the beginning of this section, at most $\Opt$.
\paragraph{Consistency.}
We can easily answer the natural question ``when the proposed estimation scheme is consistent'', meaning that for every $\epsilon\in(0,1)$, it allows to
achieve arbitrarily small $\epsilon$-risk, provided that $K$ is large enough.
Specifically, {if we denote} $G(x)=g^Tx+c$, from Proposition \ref{propwearrive} it is immediately seen that a sufficient condition for consistency is the existence of
$\bar{f}\in\bR^d$ such that
$\bar{f}^TAx=g^Tx$ for all $x\in\X-\X$, or, equivalently,  that $g$ is orthogonal to the intersection of the  kernel of $A$ with the linear span of $\X-\X$. Indeed, under this assumption, for every fixed $\epsilon\in(0,1)$
we clearly have
$\lim_{K\to\infty}\widehat{\Phi}(\bar{f})=0$, implying that $\lim_{K\to\infty}\Opt=0$, with $\widehat{\Psi}$ and $\Opt$ given by (\ref{bestprob}).
The condition in question is necessary for consistency as well,
since when the condition is violated, we have $Ax'=Ax''$ for properly selected $x',x''\in \X$ with $G(x')\neq G(x'')$, making low risk recovery of $G(x)$, $x\in\X$,  impossible already in the case of zero noise observations (i.e., those where the observation stemming from signal $x\in\X$ is identically equal to $Ax+a$)\footnote{Note that in the Gaussian case
with $M(x)$ depending on $x$ the above condition is, in general, not necessary for consistency, since a nontrivial information on $x$ (and thus on $G(x)$)
can, in principle,  be extracted from the covariance matrix $M(x)$  which can be estimated from observations.}.

\paragraph{Direct product case.}
Further simplifications are possible in the {\sl direct product case}, where, in addition to what was assumed in the beginning of Section \ref{estlinformA},
 \begin{itemize}
 \item $\E_X=\E_U\times \E_V$ and $X=U\times V$, with convex compact sets $U\subset \E_U=\bR^{n_u}$ and $V\subset E_V=\bR^{n_v}$,
 \item $\A(x=(u,v))=[Au+a,M(v)]:U\times V\to \bR^d\times \bS^d$, with $M(v)\succeq0$ for $v\in V$,
 \item $G(x=(u,v))=g^Tu+c$ depends solely on $u$, and

 \end{itemize}
 It is immediately seen that in the direct product case problem (\ref{bestprob}) reads
\begin{equation}\label{allreads}
\Opt=\min_{f\in\bR^d}\left\{\half{\left[\phi_U(A^Tf-g)+\phi_U(-A^Tf+g)\right]}+\max_{v\in V}
\left[2K^{-1}\ln(2/\epsilon)f^TM(v)f\right]^{1/2}\right\},
\end{equation}
 where
 \[
 \phi_U(h)=\max_{u\in U}u^Th.
 \]
 Assuming $\bigcap_{v\in V}\Ker(M(v))=\{0\}$, the problem is solvable, and its optimal solution $f_*$ gives rise to the
 affine estimate
 $$
 \widehat{g}_*(\omega^K)={1\over K}\sum_if_*^T\omega_i+\varkappa_*,\;\;\;\varkappa_*=\half[\phi_U(-A^Tf_*+g)-\phi_U(A^Tf_*-g)]-a^Tf_*-c,
 $$
 with $\epsilon$-risk $\leq\Opt$.
\paragraph{Near-optimality.}
In addition to the assumption that we are in the direct product case, assume  for the sake of simplicity, that $M(v)\succ0$ whenever $v\in V$. In this case (\ref{bestprob}) reads
$$
\Opt=\min_f\max_{v\in V} \left\{\Theta(f,v):=\half[\phi_U(A^Tf-g)+\phi_U(-A^Tf+g)]+
\left[2K^{-1}\ln(2/\epsilon)f^TM(v)f\right]^{1/2}\right\},
$$
whence, taking into account that $\Theta(f,v)$ clearly is convex in $f$ and concave in $v$, while $V$ is a convex compact set, by Sion-Kakutani Theorem we get also
\begin{equation}\label{sept20}
\Opt=\max_{v\in V} \left\{\Opt(v)=\min_f \half[\phi_U(A^Tf-g)+\phi_U(-A^Tf+g)]+
\left[2K^{-1}\ln(2/\epsilon)f^TM(v)f\right]^{1/2}\right\}.
\end{equation}
Now consider the problem of recovering $g^Tu$ from observation $\omega_i$, $1\leq i\leq K$, independently of each other sampled from $\N(Au+a,M(v))$, where  unknown $u$ is known to belong to $U$ and $v\in V$ is known. Let $\rho_\epsilon(v)$ be the minimax $\epsilon$-risk of the recovery:
$$
\rho_\epsilon(v)=\inf_{\widehat{g}(\cdot)}\left\{\rho:\Prob_{\omega^K\sim[\N(Au+a,M(v))]^K}\{\omega^K: |\widehat{g}(\omega^K)-g^Tu|>\rho\}\leq\epsilon\,\,\forall u\in U\right\},
$$
where $\inf$ is taken over all Borel functions $\widehat{g}(\cdot):\bR^{Kd}\to\bR$. Invoking \cite[Proposition 4.1]{JN2009}, it is immediately seen that whenever $\epsilon<{\half}$, one has
$$
\rho_\epsilon(v)\geq{{q_\N(1-\epsilon)\over \sqrt{2\ln (2/\epsilon)}}}\Opt(v)
$$
{where $q_\N(s)$ is the $s$-quantile of the standard normal distribution.}
Since the family of all sub-Gaussian, with parameters $(Au+a,M(v))$, $u\in U$, $v\in V$, distributions on $\bR^d$ contains all Gaussian distributions $\N(Au+a,M(v))$ induced by $(u,v)\in U\times V$, we arrive at the following conclusion:
\begin{proposition}\label{nearoptimality}
 In the just described situation,
the minimax optimal $\epsilon$-risk
 $$\Risk^{\hbox{\tiny\rm opt}}_{\epsilon}(K)=\inf\limits_{\widehat{g}(\cdot)} \Risk_\epsilon(\widehat{g}(\cdot)),$$
  of recovering $g^Tu$ from $K$-repeated i.i.d. sub-Gaussian, with parameters $(Au+a,M(v))$, $(u,v)\in U\times V$,  random observations is within a moderate factor of the upper bound $\Opt$ on the $\epsilon$-risk, taken w.r.t. the same data, of the affine estimate $\widehat{g}_*(\cdot)$ yielded by an optimal solution to {\rm (\ref{allreads})}. Namely,
 $$
 \Opt\leq {{ \sqrt{2\ln (2/\epsilon)}\over q_\N(1-\epsilon)}}\Risk^{\hbox{\tiny\rm opt}}_{\epsilon}
$$
with the ``near-optimality factor'' ${\sqrt{2\ln (2/\epsilon)}\over q_\N(1-\epsilon)}\to 1$ as $\epsilon\to 0$.{\footnote{It is worth mentioning that
in a more general setting of ``good observation schemes,'' described in \cite{JN2009}, the $\epsilon$-risk  $\Opt$ of the affine estimate constructed following the rules in Section \ref{sectrepeated} satisfies the bound
\[
 \Opt\leq {2\ln(2/\epsilon)\over \ln\left({1\over4\epsilon}\right)}\Risk^{\hbox{\tiny\rm opt}}_{\epsilon}
\]
where $\Risk^{\hbox{\tiny\rm opt}}_{\epsilon}$ is the corresponding minimax risk. }}
\end{proposition}

\subsection{Numerical illustration}
{In this section we consider the problem of estimating} a linear form of signal $x$ known to belong to a given convex compact subset $X$ via indirect observations $Ax$ affected by sub-Gaussian ``relative noise.'' Specifically, our observation is
$$
\omega\sim\SG(Ax,M(x))
$$
where
\begin{equation}\label{setup1}
x \in X=\left\{x\in\bR^n:0\leq x_j\leq j^{-\alpha},1\leq j\leq n\right\},\,\,
M(x)=\sigma^2\sum_{j=1}^nx_j\Theta_j.
\end{equation}
Here $A\in\bR^{d\times n}$ and $\Theta_j\in\bS^d_+$, $j=1,...,n$, are given matrices. In other words, we are in the situation where small signal results in low observation noise. The linear form to be recovered from observation
$\omega$ is $G(x)=g^Tx$. The entities $g,A,\{\Theta_j\}_{j=1}^n$ and reals  $\alpha\geq0$ (``degree of smoothness''), $\sigma>0$ (``noise intensity'') are parameters of the estimation problem we intend to process.
Parameters $g,A,\Theta_j$ are generated as follows:
 \begin{itemize}
 \item $g\geq0$ is selected at random and then normalized to have $\max\limits_{x\in X}g^Tx=2$;
 \item we consider the case of $n>d$ (``deficient observations''); the $d$ nonzero singular values of $A$ were  set to $\theta^{-{i-1\over d-1}}$, $1\leq i\leq d,$ where ``condition number''  $\theta\geq1$ is a parameter; the orthonormal systems $U$ and $V$ of the first $d$ left and, respectively, right singular vectors of $A$ were drawn at random from rotationally invariant distributions;
 \item positive semidefinite $d\times d$ matrices $\Theta_j$ are orthogonal projectors on randomly selected subspaces in $\bR^d$ of dimension $\lfloor d/2\rfloor$;
 \item in all experiments, we deal with single-observation case $K=1$.
 \end{itemize}
 Note that $X$ possesses $\geq$-largest point $\bar{x}$, whence $M(x)\preceq M(\bar{x})$ whenever $x\in X$; as a result, sub-Gaussian distributions with matrix parameter $M(x)$, $x\in X$, can be thought also to have matrix parameter $M(\bar{x})$. One of the goals of the present experiment is to {compare the risk of the affine estimate in the above model to its performance in the ``envelope model'' $\omega\sim\SG(Ax,M(\bar{x}))$, where the fact that small signals result in low-noise observations is ignored. \par

We present in Figure \ref{v_vs_m} the results of the experiment in which for a given set of parameters $d,\, n,\, \alpha,\, \theta$ and $\sigma$ we generate 100 random estimation problems -- collections $\{g,A,\Theta_j,\,j\leq d\}$. For each problem we compute $\epsilon(=0.01)$-risks of  two affine in $\omega$ estimates of $g^Tx$ as yielded by optimal solution to (\ref{bestprob}): the first -- for the problem described above (the left boxplot in each group), and the second -- for the aforementioned ``direct product envelope'' of the problem, where the mapping $x\mapsto M(x)$ is
 replaced with $x\mapsto \widehat{M}(x):=M(\bar{x})$ (the right boxplot). Note the ``noise amplification'' effect (the risk is about 20 times the level $\sigma$ of the observation noise) and significant variability of risk across the experiments. Seemingly, both these phenomena are due to deficient observation model ($n>d$) combined with ``random interplay'' between the directions of coordinate axes in $\bR^m$ (along these directions, $X$ becomes more and more thin) and the orientation of the kernel of $A$.

\begin{figure}[h!]
\centering
\includegraphics[scale=0.7]{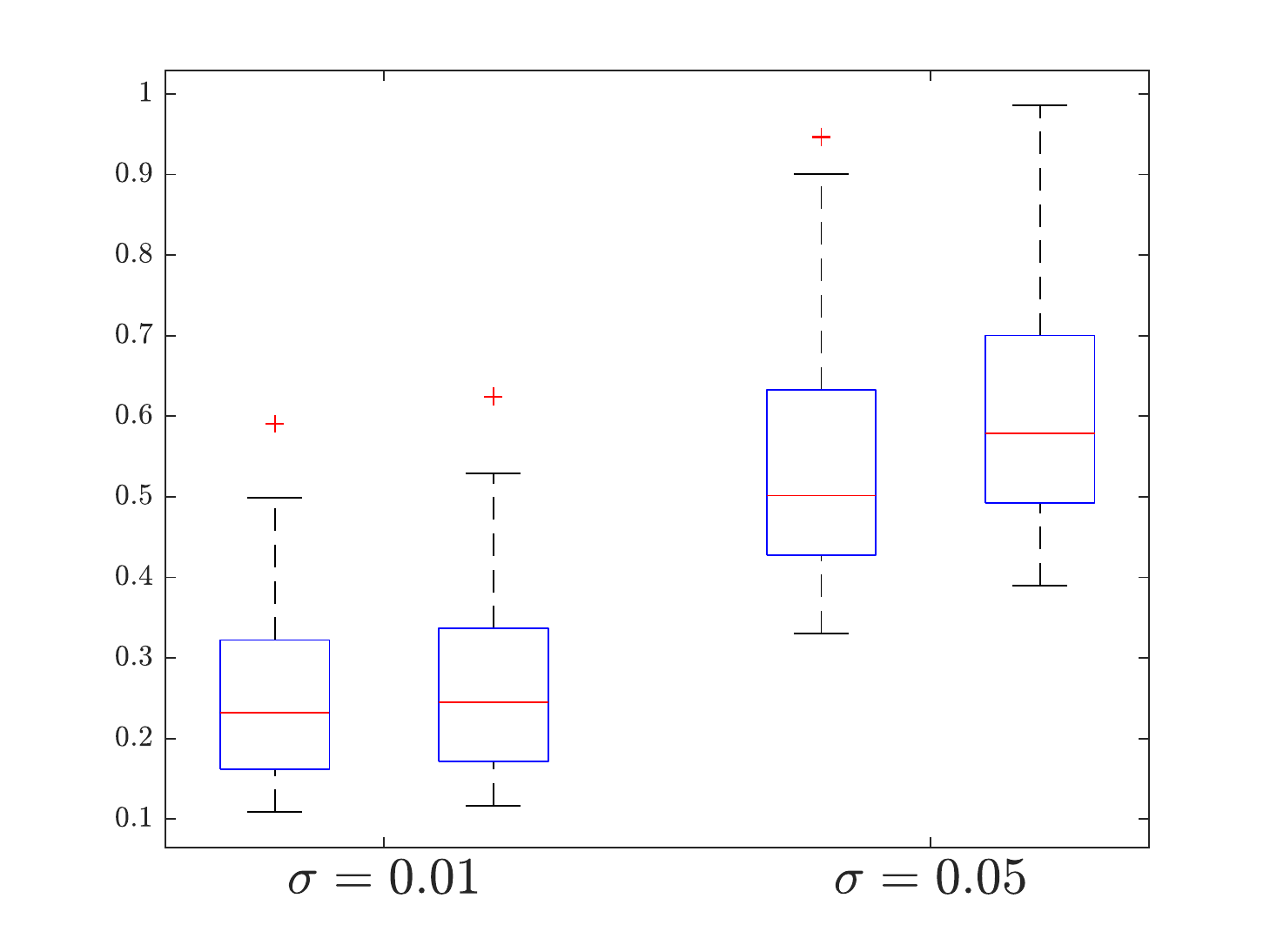}
\caption{Empirical distribution of the 0.01-risk of affine estimation over 100 estimation problems.
$\left[d=32,\ai{m}{n}=48,\alpha=2,\theta=2\right]$ for $\sigma=0.01$ and $\sigma=0.05$. In each group, distribution of risks for the problem with $\omega\sim\SG(Ax,M(x))$ on the left, for the problem with $\omega\sim\SG(Ax,M(\bar{x}))$ -- on the right. }
\label{v_vs_m}
\end{figure}
\section{Quadratic lifting and estimating quadratic forms}\label{sec:quadratic}
In this section we apply the approach in Section \ref{estlinform} to the situation where, given
an i.i.d. sample
$
\zeta^K=[\zeta_1;...;\zeta_K],\,\zeta_i\in\bR^d,
$
with distribution $P_x$ of $\zeta_i$ depending on an unknown ``signal'' $x\in X$,
 our goal is to estimate a quadratic functional $q(x)=x^TQx+c^Tx$ of the signal. We consider two situations -- the {\em Gaussian case,} where $P_{x}$ is a Gaussian distribution with parameters affinely depending on $x$, and {\em discrete case} where $P_{x}$ is a discrete distribution corresponding to the probabilistic vector $Ax$, $A$ being a given stochastic matrix.
Our estimation strategy is to apply the techniques developed in Section \ref{estlinform} to {\em quadratic liftings} $\omega$ of actual observations $\zeta$ (e.g., $\omega_i=(\zeta_i,\zeta_i\zeta_i^T)$ in the Gaussian case), so that the resulting estimates are affine functions of $\omega$'s. We first focus on implementing this program in the Gaussian case.
 \subsection{Estimating quadratic forms, Gaussian case}\label{Qlift:Gaussian}
In this section we focus on the problem  as follows. Given are
\begin{itemize}
\item  a nonempty bounded set $U\subset\bR^m$ and a nonempty convex compact set $V\subset\bR^k$,
\item an affine mapping $v\mapsto M(v):\bR^k\to\bS^d$ which maps $V$ onto convex compact subset $\V$ of $\bS^d_+$;
\item an affine mapping $u\mapsto A[u;1]:\bR^m\to\Omega=\bR^d$, where $A$ is a given $d\times(m+1)$ matrix,
\item a ``functional of interest''
\begin{equation}\label{fofint}
F(u,v)=[u;1]^TQ[u;1]+q^Tv:\;\bR^m\times\bR^k\to\bR,
\end{equation}
where $Q$ and $q$ are known $(m+1)\times (m+1)$ symmetric matrix and $k$-dimensional vector, respectively.
\item a tolerance $\epsilon\in(0,1)$.
\end{itemize}
We  observe an i.i.d. sample
$
\zeta^K=[\zeta_1;...;\zeta_K],\,\zeta_i\in\bR^d,
$
with Gaussian distribution $P_{u,v}$  of $\zeta_i$ depending on an unknown ``signal'' $(u,v)$ known to belong to $U\times V$: $P_{u,v}=\N(A[u;1],M(v))$.
 Our goal is to estimate $F(u,v)$ from
observation $\zeta^K$.
\par
The $\epsilon$-risk $\Risk_\epsilon(\widehat{g})$  of  a candidate estimate $\widehat{g}(\cdot)$ -- a Borel real-valued function on $\bR^{Kd}$ -- is defined as the smallest $\rho$ such that
$$
\forall ((u,v)\in U\times V): \Prob_{\zeta^K\sim P_{u,v}^K}\{|\widehat{g}(\zeta^K)-F(u,v)|>\rho\}\leq\epsilon.
$$
\subsubsection{Construction}
Our course of actions is as follows.
\begin{itemize}
\item
We specify  convex compact subset  $\Z\subset \bS^{m+1}$ such that
\begin{equation}\label{wesetcalZ}
\forall u\in U: [u;1][u;1]^T\in\Z\subset \Z^+=\{Z\in\bS^{m+1}_+:Z_{m+1,m+1}=1\},
\end{equation}
matrix $\Theta_*\in \bS^d$ and real $\delta\in[0,2]$ such that $\Theta_*\succ0$ and
\[
\forall \Theta\in\V: \Theta \preceq\Theta_* \mbox{and} \;\;\|\Theta^{1/2}\Theta_*^{-1/2}-I\|\leq \delta;
\]
(cf. section \ref{sect:GaussLift}).
\item We set $x(u,v)=(v, [u;1][u;1]^T)$, and $X=\{(v, [u;1][u;1]^T):\;u\in U, v\in V\}$, so that
\[
X\subset\X:=V\times \Z\subset \E_X:=\bR^K\times\bS^{m+1}.
\]
We select $\gamma\in (0,1)$ and set
\be
\H_\gamma&=&\{H\in\bS^d: -\gamma\Theta_*^{-1}\preceq H\preceq \gamma\Theta_*^{-1}\},\;\F=\bR^d\times \H_\gamma\subset \E_F=\bR^d\times \bS^d,\nn
\M&=&\V\times BZB^T\subset \E_M=\bS^d\times \bS^{d+1}, \;\;B=[A;e^T_{m+1}]
\ee{bmat}
where $e_{m+1}$ being the $(m+1)$-th canonic basis vector of $\bR^{m+1}$.
\item When adding to the above entities function $\Phi(\cdot;\cdot))$, as defined in \rf{PhisubG}, we conclude by Proposition \ref{propGausslift0} that $\M, \F$ and $\Phi(\cdot;\cdot)$ form a regular data such that for all $(u, v)\in U\times V$  and  $(h,H)\in \F$,
\be
\ln\left(
\bE_{\zeta\sim P_{u,v}}\left\{\exp\{\langle(h,H),(\zeta,\zeta\zeta^T)\rangle \}\right\}
\right)\leq \Phi\left(h,H;M(v),B[u;1][u;1]^TB^T\right)
\ee{anatoli1}
where the inner product $\langle \cdot,\cdot\rangle$  on $\E_F$ is defined as  $\langle (h,H),(g,G)\rangle=h^Tg+{1\over 2}\Tr(HG)$, so that $\langle(h,H),(\zeta,\zeta\zeta^T)\rangle=h^T\zeta+\half\zeta^TH\zeta$.
\par
Observe that $\A(x=(v,\underbrace{[u;1][u;1]^T}_{Z}))=(M(v), BZB^T)$ is an affine mapping which maps $\X$ into $\M$,
and $G(x):\E_X\to \bR$,
\[
G(x)=\Tr(QZ)+q^Tv=[u;1]^TQ[u;1]+q^Tv
\]
is a linear functional on $\E_X$.
\end{itemize}
As a result of the above steps, we get at our disposal entities $\E_X,\E_M,\E_F,\F,\M,\Phi,X,\X,\A(\cdot),G(\cdot)$ and $\epsilon$ participating in the setup described in Section \ref{linformsit}, and it is immediately seen that
these entities meet all the requirements imposed by this setup.
The bottom line is that the estimation problem stated in the beginning of this section reduces to the problem considered in Section  \ref{estlinform}.

\subsubsection{The result}\label{sec:qg-result}
When applying to the resulting data  Proposition \ref{corinfrepeated} (which is legitimate, since $\Phi$ in \rf{PhisubG} clearly satisfies (\ref{Phisatisfies})),
we arrive at the result as follows:
\begin{proposition}\label{melmelmel} In the just described situation, let us set
\begin{equation}\label{newnewnewnewnew}
\begin{array}{l}
\widehat{\Psi}_+(h,H)=\max\limits_{(v,Z)\in V\times\Z}\inf\limits_{{\alpha>0,\,\atop
-\gamma\alpha\Theta_*^{-1}\preceq H\preceq \gamma\alpha\Theta_*^{-1}}}
\left\{\alpha\Phi\left({h\over \alpha},{H\over \alpha};M(v),BZB^T\right)
-G(v,Z)+{\alpha\over K}\ln({2\over\epsilon})\right\},
\\
\widehat{\Psi}_-(h,H)=\max\limits_{(v,Z)\in V\times\Z}\inf\limits_{{\alpha>0,\atop
-\gamma\alpha\Theta_*^{-1}\preceq H\preceq \gamma\alpha\Theta_*^{-1}}}
\left\{\alpha\Phi\left(-{h\over \alpha},-{H\over \alpha};M(v),BZB^T\right)
\quad+G(v,Z)+{\alpha\over K}\ln({2\over \epsilon})\right\}.
\end{array}
\end{equation}
\noindent
so that the functions $\widehat{\Psi}_\pm(h,H):\bR^d\times\bS^d\to\bR$ are convex. Furthermore, whenever $\bar{h},\bar{H},\bar{\rho},\bar{\varkappa}$ form a feasible solution to the system of convex constraints
\begin{equation}\label{sssuchthatthat}
\begin{array}{rcl}
\widehat{\Psi}_+(h,H)\leq \rho-{\varkappa},\,\,
\widehat{\Psi}_-(h,H)\leq\rho+{\varkappa}\\
\end{array}
\end{equation}
in variables $(h,H)\in\bR^d\times \bS^d$,
$\rho\in\bR$, $\varkappa\in\bR$, setting
\begin{equation}\label{theestimate}
\widehat{g}(\zeta^K:=(\zeta_1,...,\zeta_K))={1\over K}\sum_{i=1}^K\left[h^T\zeta_i+\half\zeta_i^TH\zeta_i\right]+\bar{\varkappa},
\end{equation}
we get an estimate of the functional of interest $F(u,v)=[u;1]^TQ[u;1]+q^Tv$ via $K$ independent observations
$$
\zeta_i\sim\N(A[u;1],M(v)),\,i=1,...,K,
$$
with $\epsilon$-risk not exceeding $\bar{\rho}$:
\begin{equation}\label{resriskbound}
\forall (u,v)\in U\times V: \Prob_{\zeta^K\sim[\N(A[u;1],M(v))]^K}\left\{|F(u,v)-\widehat{g}(\zeta^K)|>\bar{\rho}\right\}\leq \epsilon.
\end{equation}
In particular,  setting for $(h,H)\in\bR^d\times\bS^d$
\begin{equation}\label{varkapparho}
\bar{\rho}=\half\left[\widehat{\Psi}_+(h,H)+\widehat{\Psi}_-(h,H)\right],\;\;\;
\bar{\varkappa}=\half \left[\widehat{\Psi}_-(h,H)-\widehat{\Psi}_+(h,H)\right],
\end{equation}
we obtain an estimate {\rm\rf{theestimate}} with  $\epsilon$-risk not exceeding $\bar{\rho}$.
\end{proposition}
For proof, see Section \ref{proof:melmelmel}.
\begin{remark}\label{remW}{\rm In the situation described in the beginning of this section, let a set $W\subset U\times V$ be given, and assume we are interested in recovering functional of interest (\ref{fofint}) at points $(u,v)\in W$ only. When reducing the ``domain of interest''  to $W$, we hopefully can reduce the  $\epsilon$-risk of recovery.  Assuming that we can point out a convex compact set $\W\subset V\times\Z$ such that
$$
(u,v)\in W\Rightarrow (v,[u;1][u;1]^T)\in \W.
$$
it can be straightforwardly verified that in this case   the conclusion of Proposition \ref{melmelmel} remains valid when the set $V\times\Z$ in (\ref{newnewnewnewnew})
is replaced with $\W$, and the set $U\times V$ in (\ref{resriskbound}) is replaced with $W$.
 This modification enlarges the feasible set of
(\ref{sssuchthatthat}) and thus reduces the attainable risk bound.}
\end{remark}
\paragraph{Discussion.}
When estimating quadratic forms from $K$-repeated observations $\zeta^K=[\zeta_1;...;\zeta_K]$ with i.i.d. $\zeta_i$  
we applied ``literally'' the construction of Section \ref{sectrepeated}, thus restricting ourselves with estimates affine in quadratic liftings $\omega_i=(\zeta_i,\zeta_i\zeta_i^T)$ of $\zeta_i$'s. As an alternative to such ``basic'' approach, let us consider estimates which are affine in the ``full'' quadratic lifting $\omega=(\zeta^K,\zeta^K[\zeta^K]^T)$ of $\zeta^K$,  thus extending the family of candidate estimates (what is affine in $\omega_1,...,\omega_K$, is affine in $\omega$, but not vice versa, unless $K=1$). Note that this alternative is covered by our approach -- all we need, is to replace the original components  $d$, $M(\cdot)$, $\V$, $A$ of the setup of this section with their extensions
$$\begin{array}{c}d^+=Kd,\;M^+(v)=\Diag\{\underbrace{M(v),...,M(v)}_{K}\},\\
\V^+=M^+(V)=\{\Theta=\Diag\{M(v),...,M(v)\},v\in V\},\;A^+=[A;...;A],
\end{array}
$$
 and set $K$ to 1.
\par
It is easily seen that {such modification} can only reduce the risk of the resulting estimates, the price being the increase in design dimension (and thus in computational complexity) of the optimization problems yielding the estimates.
To illustrate the difference between two approaches, consider the situation (to be revisited in Section \ref{sec:direct-sim}) where we are interested to recover the energy $u^Tu$ of a signal $u\in\bR^m$ from observation
\begin{equation}\label{gaussobs}
\zeta=u+\xi,\;\;\xi\sim\N(0,\Theta)
\end{equation}
where $\Theta$ is (unknown) diagonal matrix with diagonal entries from the range $[0,\sigma^2]$, and a priori information about $u$ is that $\|u\|_2\leq R$ for some known $R$. Assume that (cf. Section \ref{directobsproc}) $m\geq16\ln(2/\epsilon)$, where $\epsilon\in(0,1)$ is a given reliability tolerance and that $R^2\geq m\sigma^2$.}
Under these assumptions one can easily verify that in the single-observation case the $\epsilon$-risks of both the ``plug-in'' estimate $\zeta^T\zeta$ and of the estimate yielded by the proposed approach are, up to absolute constant factors, the same as the  optimal $\epsilon$-risk, namely, $O(1)\R$, $\R=\sigma^2m+\sigma R\sqrt{\ln(2/\epsilon)}$.
Now let us look at the case $K=2$ where we observe two independent copies, $\zeta_1$ and $\zeta_2$, of observation (\ref{gaussobs}). Here the $\epsilon$-risks of the ``naive'' plug-in estimate ${1\over 2}[\zeta_i^T\zeta_1+\zeta_2^T\zeta_2]$,
and of the estimate obtained by applying our ``basic'' approach with $K=2$ are just by absolute constant factors better than in the single-observation case -- both these risks still are $O(1)\R$.
In contrast to this, an ``intelligent'' plug-in 2-observation estimate $\zeta_1^T\zeta_2$ has risk $O(1)\sigma (R+\sigma\sqrt{m})\sqrt{\ln(2/\epsilon)}$ whenever $R\geq0$, which is much smaller than $\R$ when $m\gg \ln(2/\sigma)$ and  $R\sqrt{\ln(2/\epsilon)}\ll \sigma m$.
It is easily seen that with the outlined alternative implementation, our approach also results in estimate with ``correct'' $\epsilon$-risk $O(1)\sigma (R+\sigma\sqrt{m})\sqrt{\ln(2/\epsilon)}$.

\subsubsection{Consistency}
We are about to present a simple sufficient condition for the estimator suggested by Proposition \ref{melmelmel} to be consistent, in the sense of Section {\ref{estlinformA}}. Specifically,
assume that
\begin{itemize}
\item[{A.1}.] $V=\{\bar{v}\}$ is a singleton such that $M(\bar{v})\succ0$, which allows to satisfy (\ref{56delta}) with $\Theta_*=M(\bar{v})$ and $\delta=0$, same as allows to assume w.l.o.g. that
$$
F(u,v)=[u;1]^TQ[u;1],\,\,G(x=(v,Z))=\Tr(QZ);
$$
\item[{A.2}.] the first $m$ columns of the $d\times (m+1)$ matrix $A$ are linearly independent.
\end{itemize}
The consistency of our estimation procedure is given by the following simple statement:
\begin{proposition}\label{propconsistency} In the just described situation and under assumptions A.1--2, given $\epsilon\in(0,1)$, consider the estimate
$$
\widehat{g}_{K}(\zeta^K) = {1\over K}\sum_{i=1}^K[\bar{h}^T\zeta_i+\half\zeta_i^T\bar{H}\zeta_i]+\varkappa_{K},
$$
where
$$
\varkappa_{K}=\half \left[\widehat{\Psi}_-(\bar{h},\bar{H})-\widehat{\Psi}_+(\bar{h},\bar{H})\right]
$$
and $\widehat{\Psi}_\pm=\widehat{\Psi}_\pm^{K}$ are given by {\rm (\ref{newnewnewnewnew})}. Then the $\epsilon$-risk of $\widehat{g}_{K,\epsilon}(\cdot)$ goes to 0 as $K\to\infty$.
\end{proposition}

For proof, see Section \ref{NEwProof}.

\subsection{Numerical illustration, direct observations}\label{sec:direct-sim}
\subsubsection{The problem}
Our first illustration is deliberately selected to be extremely simple: given direct noisy observation
$$
\zeta=u+\xi
$$
of unknown signal $u\in\bR^m$ known to belong to a given set $U$, we want to recover the ``energy'' $u^Tu$ of $u$; what we are interested in, is the quadratic in $\zeta$ estimate with as small $\epsilon$-risk on $U$ as possible; here $\epsilon\in(0,1)$ is a given design parameter. Note that we are in the situation where the dimension $d$ of the observation is equal to the dimension $m$ of the signal underlying observation. The details of our setup are as follows:
\begin{itemize}
\item $U$ is the ``spherical layer'' $U=\{u\in\bR^m: r^2\leq u^Tu\le R^2\}$, where $r,R$, $0\leq r<R<\infty$ are given. As a result, the ``main ingredient'' of constructions in Section \ref{Qlift:Gaussian}
 -- the convex compact subset $\Z$ of $\Z^+$ containing all matrices $[u;1][u;1]^T$, $u\in U$, see (\ref{wesetcalZ}), can be specified as
    $$
    \Z=\big\{Z\in\bS^{m+1}_+:Z_{m+1,m+1}=1,1+r^2\leq \Tr(Z)\leq 1+R^2\big\};
    $$
\item $\xi\sim \N(0,\Theta)$, with matrix $\Theta$ known to be diagonal with diagonal entries
satisfying $\theta\sigma^2\leq \Theta_{ii}\leq\sigma^2$, $1\leq i\leq d=m$, with known $\theta\in[0,1]$ and $\sigma^2>0$. In terms the setup of Section \ref{Qlift:Gaussian}, we are in the case where $V=\{v\in\bR^m: \theta\sigma^2\leq v_i\leq \sigma^2,i\leq m\}$, and $M(v)=\Diag\{v_1,...,v_m\}$;
\item the functional of interest is
$
F(u,v)=u^Tu,
$
i.e., is given by (\ref{fofint}) with $Q=I_m$ and $q=0$.
\end{itemize}
\subsubsection{Processing the problem}\label{directobsproc}
It is easily seen that in the situation in question the construction in Section \ref{Qlift:Gaussian} boils down to the following:
\begin{enumerate}
\item We lose nothing when restricting ourselves with estimates of the form
\begin{equation}\label{estimateyielded}
\widehat{g}(\zeta)={\eta\over 2}\zeta^T\zeta +\varkappa,
\end{equation}
with properly selected scalars $\eta$ and $\varkappa$;
\item 
$\eta$ and $\varkappa$ are supplied by the convex optimization problem (with just 3 variables $\alpha_{+},\alpha_{-},\eta$)
\be
\min\limits_{\alpha_{\pm},\eta}\left\{\widehat{\Psi}(\alpha_+,\alpha_-,\eta)=\half\left[\widehat{\Psi}_{+}(\alpha_{+},\eta)+\widehat{\Psi}_{-}(\alpha_{-},\eta)\right]:
\sigma^2|\eta|<\alpha_\pm\right\},
\ee{qwerty1}
where
\[
\begin{array}{rcl}
\widehat{\Psi}_+(\alpha_{+},\eta)&=&-{m\alpha_+\over 2}\ln(1-\sigma^2\eta/\alpha_+)+{m\over 2}\sigma^2(1-\theta)\max[-\eta,0]+
{m\delta(2+\delta)\sigma^4\eta^2
\over2(\alpha_+-\sigma^2|\eta|)}\\
&&+\max\limits_{r^2\leq t\leq R^2}\left[
\left[{\alpha_+\eta\over 2(\alpha_+-\sigma^2\eta)}-1\right]t\right]+\alpha_+\ln(2/\epsilon)\\
\widehat{\Psi}_-(\alpha_{-},\eta)&=&-{m\alpha_-\over 2}\ln(1+\sigma^2\eta/\alpha_-)+{m\over 2}\sigma^2(1-\theta)\max[\eta,0]+
{m\delta(2+\delta)\sigma^4\eta^2
\over2(\alpha_--\sigma^2|\eta|)}\\
&&+\max\limits_{r^2\leq t\leq R^2}\left[
\left[-{\alpha_-\eta\over 2(\alpha_-+\sigma^2\eta)}+1\right]t\right]+\alpha_-\ln(2/\epsilon),
\end{array}
\]
with $\delta=1-\sqrt{\theta}$. Specifically, the $\eta$-component of a feasible solution to \rf{qwerty1} augmented by the quantity
$$
\varkappa={1\over 2}\left[\widehat{\Psi}_-(\alpha_-,\eta)-\widehat{\Psi}_+(\alpha_+,\eta)\right]
$$
yields estimate (\ref{estimateyielded}) with $\epsilon$-risk on $U$ not exceeding $\widehat{\Psi}(\alpha_+,\alpha_-,\eta)$;

\end{enumerate}
The  ``energy estimation'' problem {where} $\xi\sim\N(0,\sigma^2I_m)$ with  $\sigma^2$ known to belong to a given range is well studied in the literature. Available results investigate analytically the interplay between the dimension $m$ of signal, the range of noise intensity $\sigma^2$ and the parameters $R,r,\epsilon$ and offer provably optimal, up to absolute constant factors, estimates. For example, consider the
case with  $r=0$, and
$\theta=0$, and assume for the sake of definiteness that $R^2\geq \sigma^2m$ (otherwise already the trivial -- identically zero -- estimate is near optimal) and that we are in ``high dimensional regime,'' i.e., $m\geq 16\ln(2/\epsilon)$.  It is well known that in this case the optimal $\epsilon$-risk, up to absolute constant factor, is $\sigma^2m+\sigma R\sqrt{\ln(2/\epsilon)}$ and is achieved, again, up to absolute constant factor, at the ``plug-in'' estimate
$\widehat{x}(\zeta)=\zeta^T\zeta
$. It is easily seen that under the circumstances similar risk bound holds true for the estimate (\ref{estimateyielded}) yielded by the optimal solution to (\ref{qwerty1}). \par
A nice property of the proposed approach is that (\ref{qwerty1}) automatically takes care of the parameters and results in estimates with seemingly near-optimal performance, as is witnessed by the numerical results we present below.
\subsubsection{Numerical results}
In the experiments we are reporting {on,} we compute, for different sets of parameters $m,\,r,\,R$ and $\theta$ ($\sigma=1$ in all experiments)  the 0.01-risk attainable by the proposed estimators in the Gaussian case -- the optimal values of the problem \rf{qwerty1}, along with ``suboptimality ratios'' of such risks to the lower bounds on the best possible under circumstances 0.01-risks.

To compute these lower bounds we use the following construction. Consider the problem of estimating $\|u\|^2_2$, $u\in U=\{u:\,r\le \|u\|_2\leq R\}$ given observation $\omega=\N(u,\vartheta I_m)$, with $\vartheta\in [\theta,1]$. Same as in Section \ref{estlinformA}, the optimal $\epsilon$-risk $\Risk^{\hbox{\tiny\rm opt}}_\epsilon$ for this problem is defined as the infimum of the $\epsilon$-risk over all estimates. Now let us select somehow the  $r_1,r_2$, $r\leq r_1<r_2\leq R$, and $\sigma_1,\sigma_2$,  $\theta\leq \sigma_1,\sigma_2\leq 1$, and let $P_1$ and $P_2$ be
two distributions of observations as follows: $P_\chi$ is the distribution of random vector $\omega=\eta+\xi$, where $\eta$  and $\xi$ are independent, $\eta$ is uniformly distributed over the sphere $\|\eta\|_2=r_\chi$, and $\xi\sim\N(0,\sigma_\chi^2I_m)$, $\chi=1,2$. It is immediately seen that if there is no test which can decide on the hypotheses $H_1: \,\omega\sim P_1$, and $H_2:\,\omega\sim P_2$ via observation $\omega$  with total risk $\leq2\epsilon$ (defined as the sum, over our two hypotheses, of probabilities to reject the hypothesis when it is true), the quantity ${r_2^2-r_1^2\over 2}$ is a lower bound on the optimal $\epsilon$-risk $\Risk^{\hbox{\tiny\rm opt}}_\epsilon$. In other words, denoting by $p_\chi(\cdot)$ the density of $P_\chi$, we have
$$
0.02<\int_{\bR^d}\min[p_1(\omega),p_2(\omega)]d\omega \Rightarrow \Risk^{\hbox{\tiny\rm opt}}_{0.01}\geq {r_2^2-r_1^2\over 2}.
$$
Now, the densities $p_\chi$ are spherically symmetric, whence, denoting by $q_\chi(\cdot)$ the univariate density of the energy $\omega^T\omega$ of observation  $\omega\sim P_\chi$, we have
$$
\int_{\bR^d}\min[p_1(\omega),p_2(\omega)]d\omega =\int_0^\infty\min[q_1(s),q_2(s)]ds,
$$
and we conclude that
\begin{equation}\label{shaovn}
0.02 <\int_0^\infty\min[q_1(s),q_2(s)]ds \Rightarrow \Risk^{\hbox{\tiny\rm opt}}_{0.01}\geq {r_2^2-r_1^2\over 2}.
\end{equation}
On a closest inspection, $q_\chi$ is the convolution of two univariate densities representable by explicit computation-friendly formulas, implying that given $r_1,r_2,\sigma_1,\sigma_2$, we can check numerically whether the premise in (\ref{shaovn}) indeed takes place; whenever this is the case, the quantity $ {r_2^2-r_1^2\over 2}$ is a lower bound on $\Risk^{\hbox{\tiny\rm opt}}_{0.01}$. In our experiments, we used a simple search strategy (not described here) aimed at crude maximizing this bound in $r_1,r_2,\sigma_1,\sigma_2$ and used the resulting lower bounds on $\Risk^{\hbox{\tiny\rm opt}}_{0.01}$ to compute the suboptimality ratios.\footnote{The reader should not be surprised by the ``singular numerical spectrum'' of optimality ratios: our lower bounding scheme was restricted to identify actual optimality ratios among the candidate values $1.05^i$, $i=1,2,...$}
\par
\begin{figure}[h!]
\begin{tabular}{cc}
\includegraphics[scale=0.55]{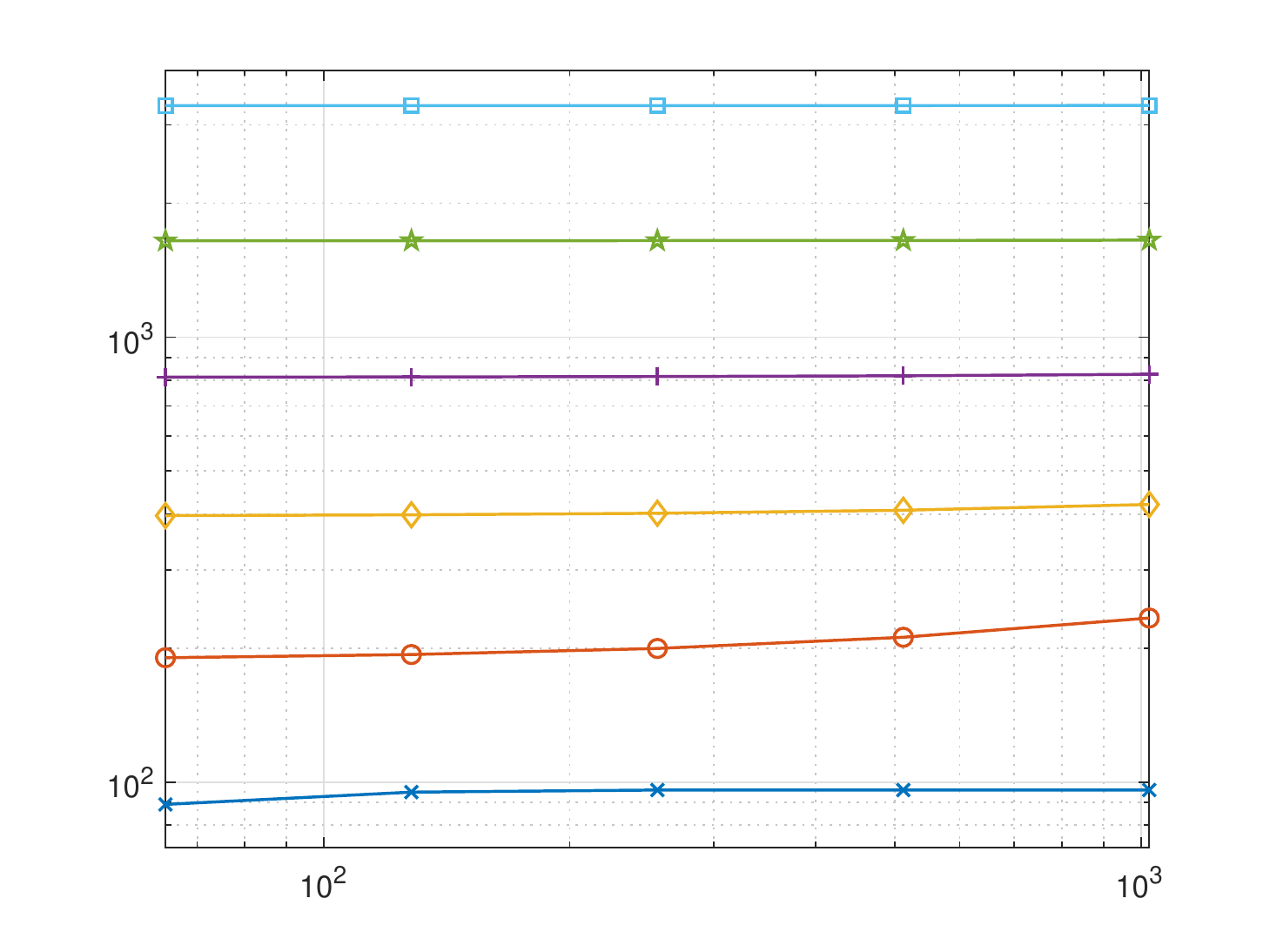}&\includegraphics[scale=0.55]{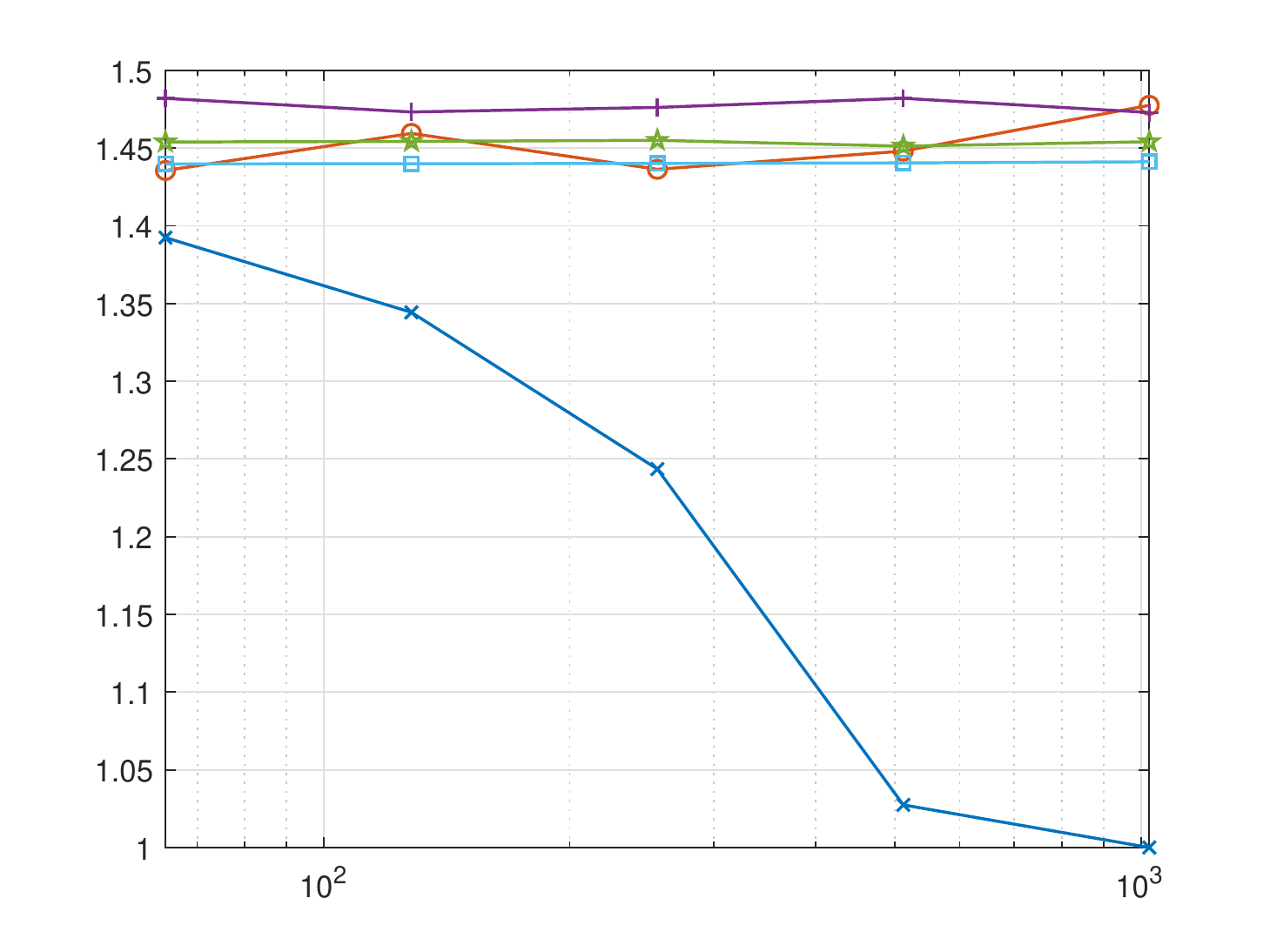}
\end{tabular}
\caption{\small Estimation risks as functions of problem dimension $m$ and $R\in\{16, 32, 64, 128, 256, 512\}$ (different curves); other parameters: $r/R=0.5$, and $\theta=1.0$.
Left plot: estimation risks; right plot: suboptimality ratios.}
\label{fig:quad1}
\end{figure}
\begin{figure}[h!]
\begin{tabular}{cc}
\includegraphics[scale=0.55]{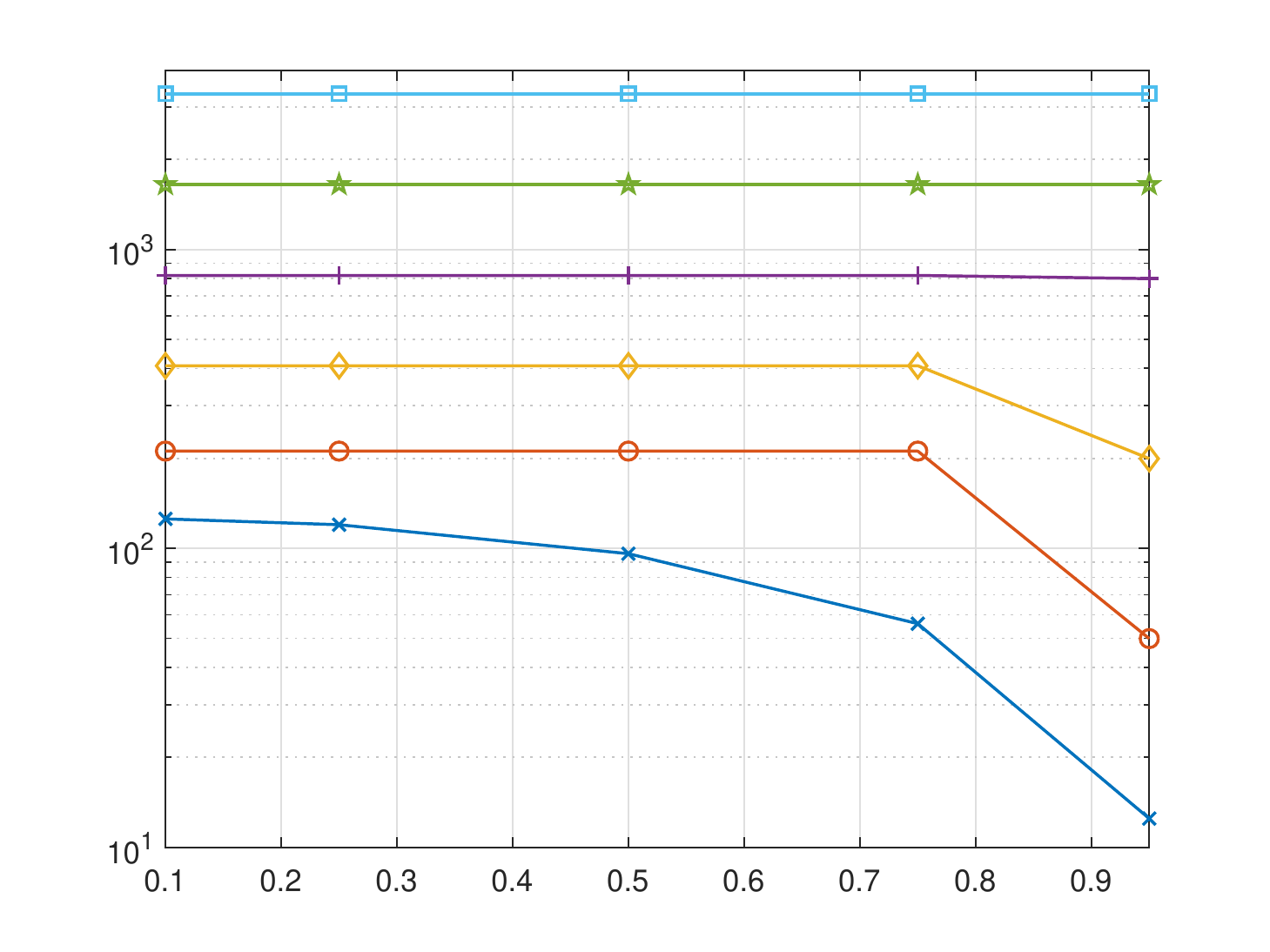}&\includegraphics[scale=0.55]{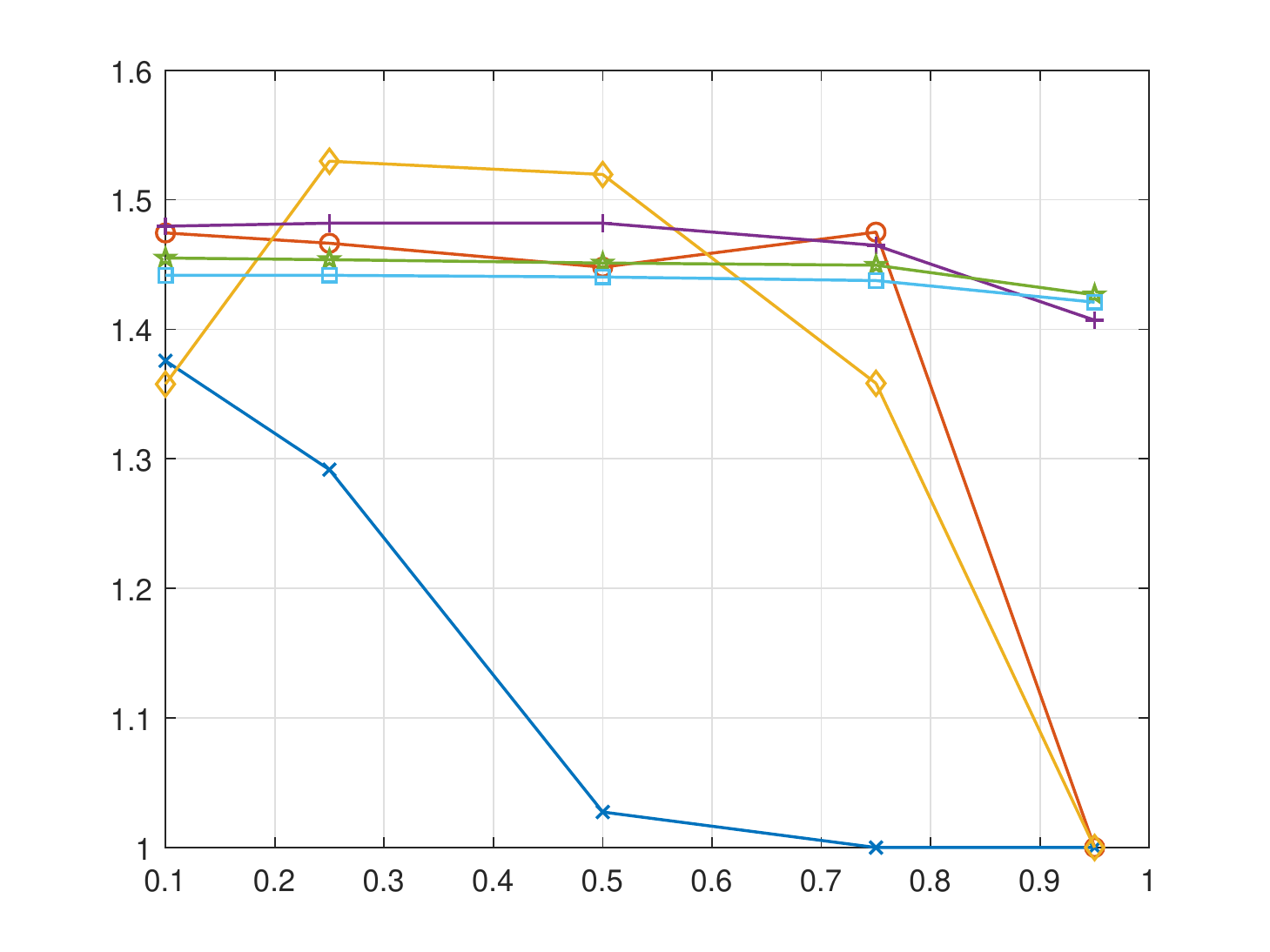}
\end{tabular}
\caption{\small Estimation risks as functions of the ratio $r/R$ and $R\in\{16, 32, 64, 128, 256, 512\}$ (different curves); other parameters: $m=512$, and $\theta=1.0$.
Left plot: estimation risks; right plot: suboptimality ratios.}
\label{fig:quad2}
\end{figure}
\begin{figure}[h!]
\begin{tabular}{cc}
\includegraphics[scale=0.55]{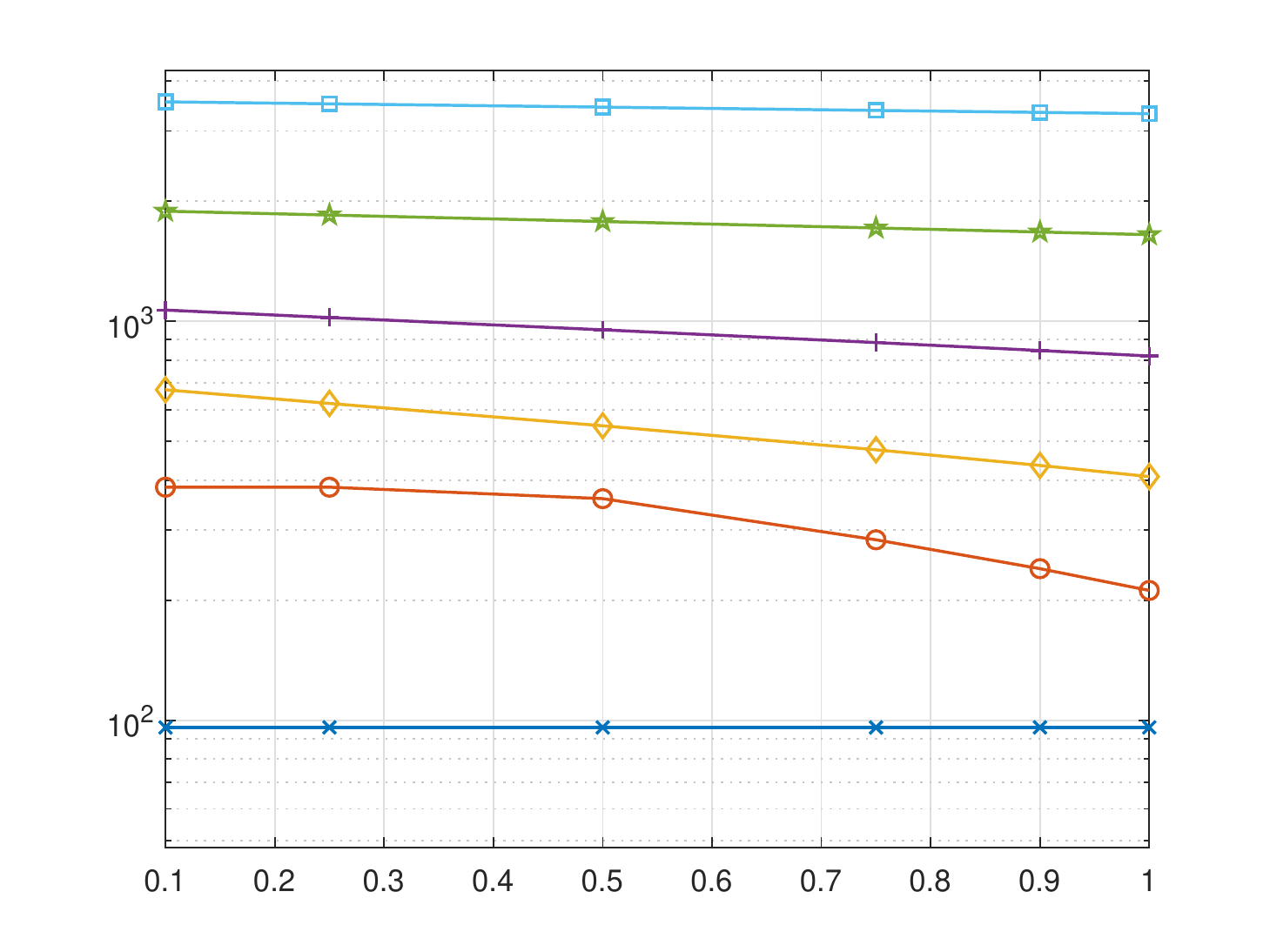}&\includegraphics[scale=0.55]{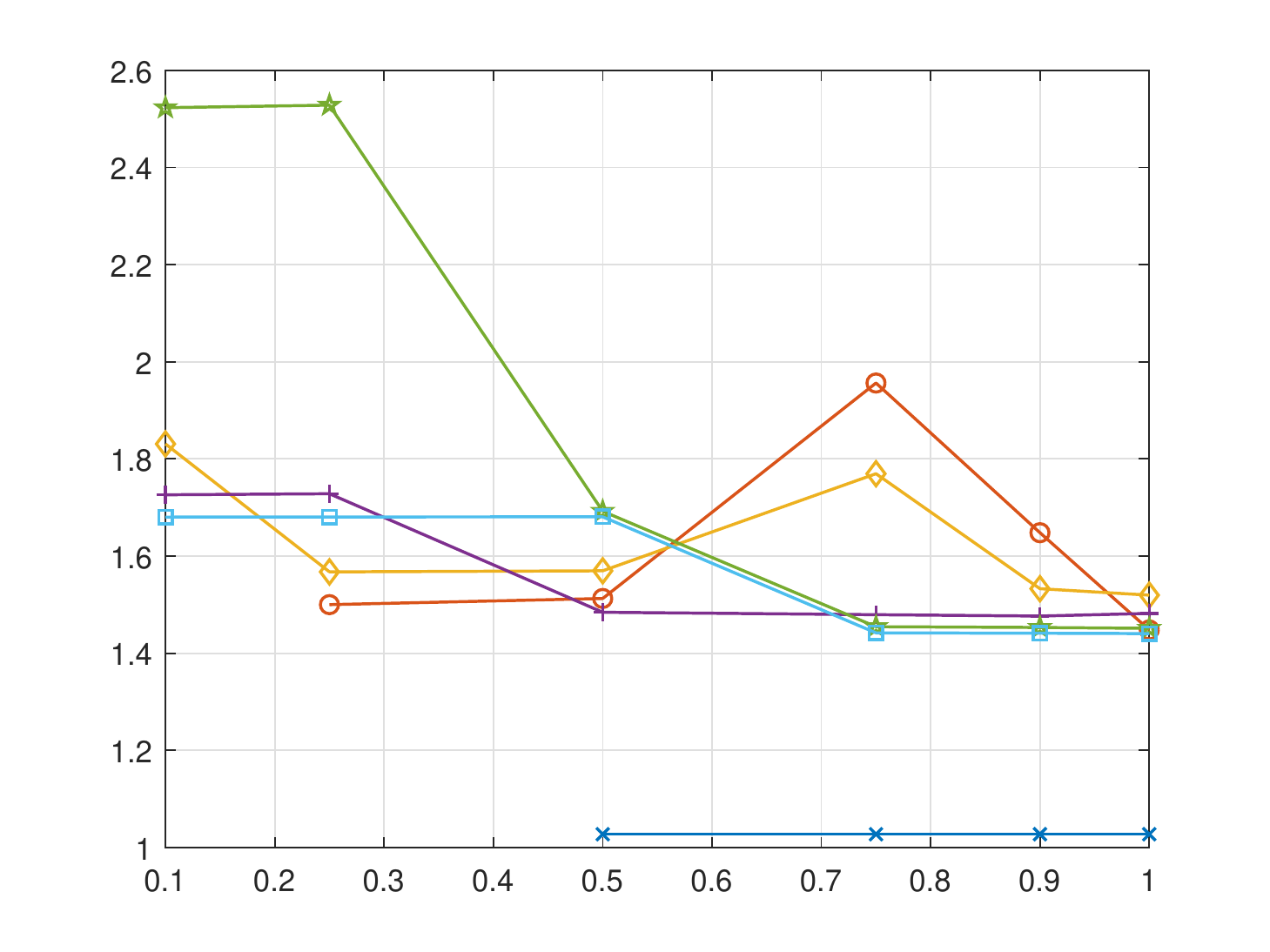}
\end{tabular}
\caption{\small Estimation risks as functions of $\theta$ and $R\in\{16, 32, 64, 128, 256, 512\}$ (different curves); other parameters: $m=512$, and $r/R=0.5$.
Left plot: estimation risks; right plot: suboptimality ratios.}
\label{fig:quad3}
\end{figure}

In Figures \ref{fig:quad1}--\ref{fig:quad3} we present some typical simulation results illustrating  dependence of risks on problem dimension $m$ (Figure \ref{fig:quad1}), on ratio $r/R$ (Figure \ref{fig:quad2}), and on parameter $\theta$ (Figure \ref{fig:quad3}). Different curves in each plot correspond to different values of the parameter $R$ varying in $\{16, 32, 64, 128, 256, 512\}$, other parameters being fixed. We believe that quite moderate values of the optimality ratios presented in the figures (these results are typical for a much larger series of experiments we have conducted) attest a rather good performance of the proposed apparatus.

\subsection{Numerical illustration, indirect observations}
\subsubsection{The problem}
The estimation problem we {address in this section} is as follows. Our observations are
\begin{equation}\label{illeq1}
\zeta=Pu+\xi,
\end{equation}
where
\begin{itemize}
\item $P$ is a given $d\times m$ matrix, with $m>d$ (``under-determined observations''),
\item $u\in\bR^m$ is a signal known to belong to a given compact set $U$,
\item $\xi\sim\N(0,\Theta)$  is the observation noise; $
\Theta$ is  positive semidefinite $d\times d$ matrix known to belong to a given convex compact set $\V\subset\bS^d_+$.
\end{itemize}
Our goal is to estimate the energy
$$
F(u)={\|u\|_2^2\over m}
$$
of the signal given a single observation (\ref{illeq1}).
\par
In our experiment, the data is specified as follows:
\begin{enumerate}
\item We assume that $u\in\bR^m$ is a discretization of a smooth function $x(t)$ of continuous argument $t\in[0;1]$: $u_i=x({i
\over m})$, $1\leq i\leq m$,  and use in the role of $U$  ellipsoid $\{u\in\bR^m: \|Su\|_2^2\leq 1\}$ with $S$ selected to make $U$ a natural discrete-time version of
the Sobolev-type ball $
\{x:[x(0)]^2+[x'(0)]^2+
\int_0^1[x''(t)]^2dt \leq 1\}$.
\item $d\times m$ matrix $P$ is of the form $UDV^T$, where $U$ and $V$ are randomly selected $d\times d$ and $m\times m$ orthogonal matrices, and the $d$ diagonal entries in diagonal $d\times m$ matrix $D$ are of the form $\vartheta^{-{i-1\over d-1}}$, $1\leq i\leq d$; the ``condition number'' ${\vartheta}$ of $P$ is a design parameter.
\item The set $\V$ of allowed values of the ``covariance'' matrices $\Theta$ is the set of all diagonal $d\times d$ matrices with diagonal entries varying in  $[0, \sigma^2]$, with the ``noise intensity'' $\sigma$ being a design parameter.
\end{enumerate}
\subsubsection{Processing the problem}
Our estimating problem clearly is covered by the setups considered in Section \ref{Qlift:Gaussian}. In terms of these setups, we specify $\Theta_*$ as $\sigma^2I_d$, $V$ as $\V$, and
$M(v)$ as the identity mapping of $\bS^d$ onto itself; the mapping $u\mapsto A[u;1]$ becomes the mapping $u\mapsto Pu$, while the set $\Z$ (which should be a convex compact subset of the set $\{Z\in\bS^{d+1}_+:Z_{d+1,d+1}=1\}$ containing all matrices of the form $[u;1][u;1]^T$, $u\in U$) becomes the set
$$
\Z=\{Z\in\bS^{d+1}_+:Z_{d+1,d+1}=1,\Tr\left(Z\Diag\{S^TS,0\}\right) \leq m\}.
$$
As suggested by Proposition \ref{melmelmel},  linear in ``lifted observation'' $\omega=(\zeta,\zeta\zeta^T)$ estimates of $F(u)={1\over m}\|u\|_2^2$ stem from the optimal solution $(h_*,H_*)$ to the convex optimization problem
\begin{equation}\label{illprob}
\Opt=\min_{h,H}\half[\widehat{\Psi}_+(h,H)+\widehat{\Psi}_-(h,H)],
\end{equation}
with $\widehat{\Psi}_{\pm}(\cdot)$ given by (\ref{newnewnewnewnew}) as applied with  $K=1$. The resulting estimate is
\begin{equation}
\label{estres}
\zeta\mapsto h_*^T\zeta+\half\zeta^TH_*\zeta +\varkappa,\,\,\varkappa=\half [\widehat{\Psi}_-(h_*,H_*)-\widehat{\Psi}_+(h_*,H_*)]
\end{equation}
and the $\epsilon$-risk of the estimate is (upper-bounded by) $\Opt$.
\par
Problem (\ref{illprob}) is a well-structured convex-concave saddle point problem and as such is beyond the ``immediate scope'' of the standard Convex Programming software toolboxes primarily aimed at solving well-structured
convex minimization problems. However, applying conic duality, one can easily eliminate in (\ref{newnewnewnewnew}) the inner maxima over $v,Z$ ro arrive at the reformulation which can be solved numerically by {\tt  CVX} \cite{cvx2014}, and this is how (\ref{illprob}) was processed in our experiments.
\subsubsection{Numerical results}
 To quantify the performance of the proposed approach, we present, along with the upper risk bounds, simple lower bounds on the best $
\epsilon$-risk achievable under the circumstances. The origin of these lower bounds is as follows. Let $w\in U$ with $t(w)=\|Pw\|_2$, and let $\rho=2\sigma q_\N(1-\epsilon)$  where $q_\N(\cdot)$ is the standard normal quantile:
\[
\Prob_{\xi\sim\N(0,1)}\{\xi\leq q_\N(p)\}=p\;\;\forall p\in(0,1).
 \]
Then for $\theta(w)=\max[1-\rho/t(w),0]$, we have $w':=\theta(w)w\in U$, and $\|Pw-Pw'\|_2\leq \rho$. The latter, due to the origin of $\rho$, implies  that there is no test which decides on the hypotheses $u=w$ and $u=w'$ via observation $Pu+\xi$, $\xi\sim\N(0,\sigma^2I_d)$, with risk $<\epsilon$.
As an immediate consequence, the quantity \[
\phi(w):={1\over 2}[\|w\|_2^2-\|w'\|_2^2]=\|w\|_2^2[1-\theta^2(w)]/2\]
is a lower bound on the $\epsilon$-risk, on $U$, of a whatever estimate of $\|u\|_2^2$. We can now try to maximize the resulting lower risk bound over $U$, thus arriving at the lower bound
 $$
 \hbox{LwBnd}=\max_{w\in U}\left\{\half\|w\|_2^2(1-\theta^2(w))\right\}.
 $$
 On a closest inspection, the latter problem is not a convex one, which does not prevent us from building its suboptimal solution.
 \par
\begin{figure}[h!]
\begin{tabular}{cc}
\includegraphics[scale=0.55]{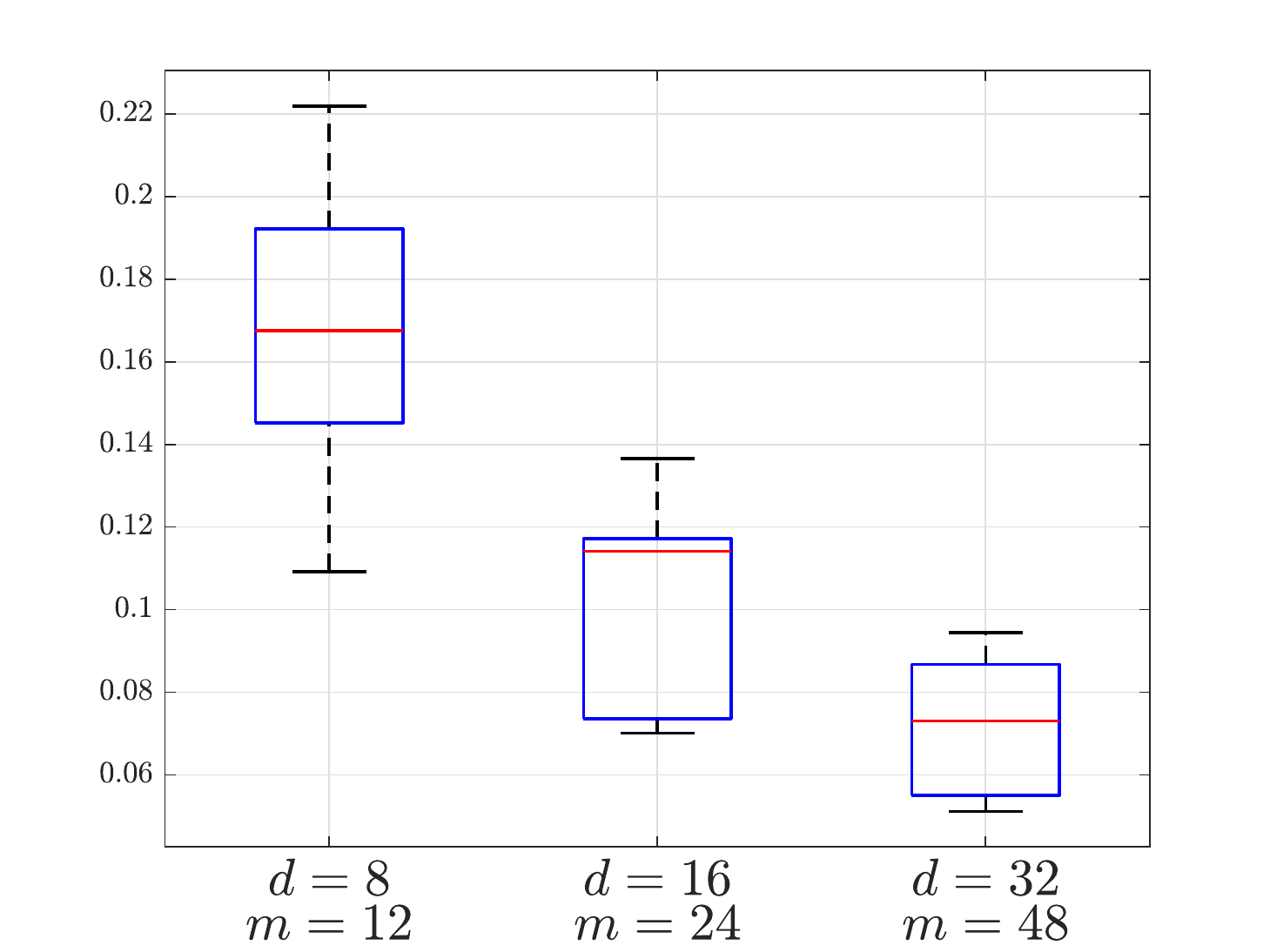}&\includegraphics[scale=0.55]{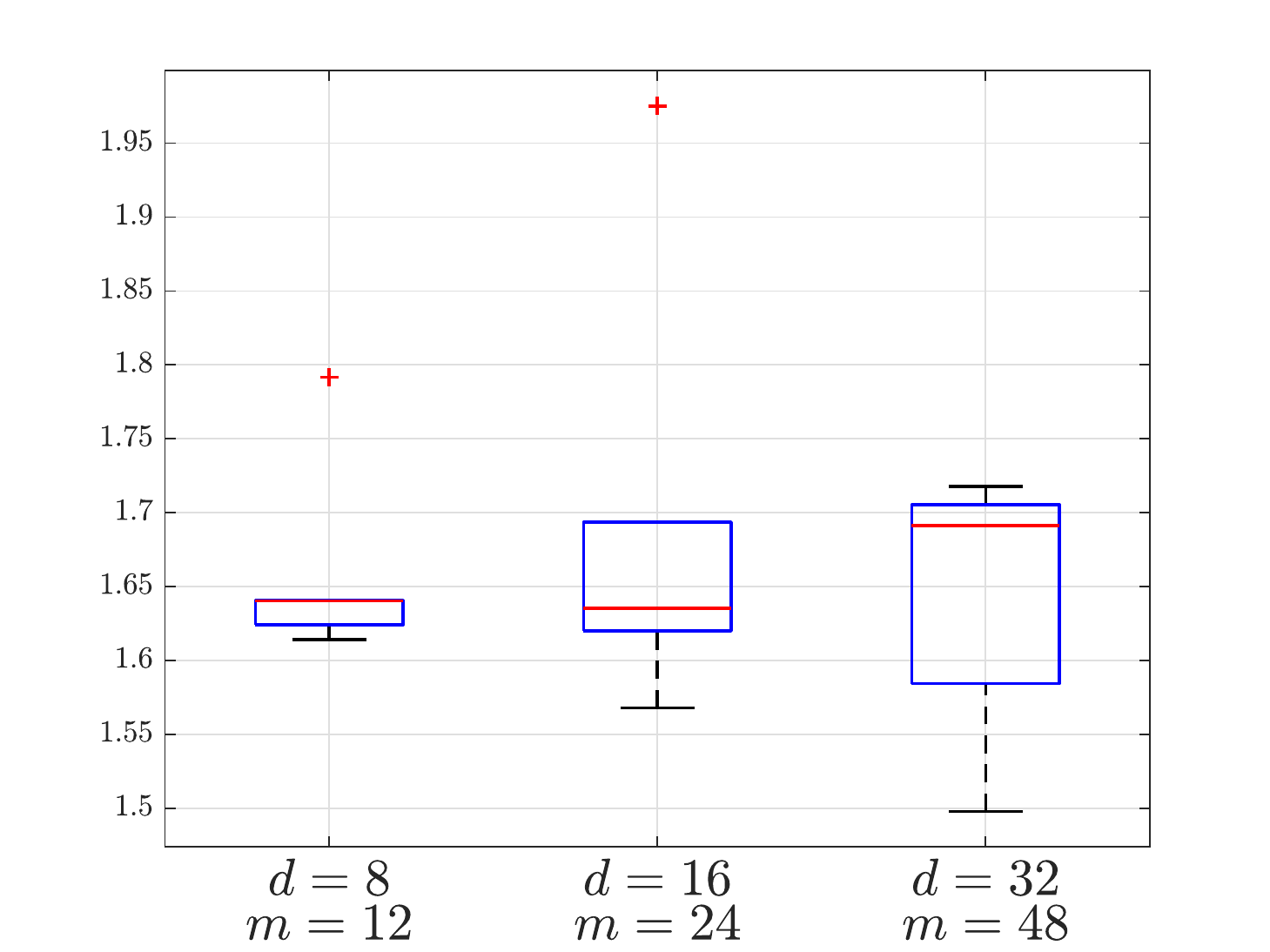}\\
$(a)$&$(b)$
\end{tabular}
\caption{\small Empirical distribution of the 0.01-risk over 20 random estimation problems, $\sigma=0.025$. $(a)$: upper risk bound $\Opt$ as in \rf{illprob}; $(b)$ corresponding suboptimality ratios. }
\label{fig:gind1}
\end{figure}

 Note that in our experiments even with fixed design parameters $d,m,\theta,\sigma$, we still deal with families of estimation problems differing from each other by their  ``sensing matrices'' $P$;  orientation of the system of right singular vectors of $P$ with respect to the axes of $U$  is random, so that these matrices varies essentially from simulation to simulation, which affects significantly the attainable estimation risks.
 We display in Figure \ref{fig:gind1} typical results of our experiments.
 We see that the (theoretical upper bounds on the) $\epsilon$-risks of our estimates, while varying significantly with the parameters of the experiment, all the time stay within a moderate factor from the lower risk bounds.

\subsection{Estimation of quadratic functionals of a discrete distribution}
In this section we consider the situation as follows: we are given an $d\times m$ ``sensing matrix'' $A$ which is stochastic -- with columns belonging to the probabilistic simplex $\Delta_d=\{v\in\bR^d:\;v\geq0,\sum_iv_i=1\}$, and a nonempty closed subset $U$ of $\Delta_m$, along with a $K$-repeated observation $\zeta^{K}=(\zeta_1,...,\zeta_{K})$ with $\zeta_i$, $1\leq i\leq K$, drawn independently across $i$ from the discrete distribution $\mu=Au_*$, where $u_*$ is an unknown probabilistic vector (``signal'') known to belong to $U$. {\sl We always assume that $K\geq 2$}. We treat a discrete distribution on $d$-point set as a distribution $P_\mu$ on the $d$ vertices $e_1,...,e_d$ of $\Delta_d$, so that possible values of $\zeta_i$ are basic orths $e_1,...,e_d$ in $\bR^d$ with $\Prob_{\zeta\sim \mu}(\zeta=e_{j})=\mu_j$. Our goal is
to recover from observation $\zeta^K$ the value at $u_*$ of a given quadratic form
$$
F(u)=u^TQu+2q^Tu.
$$
\subsubsection{Construction}
Observe that for $u\in\Delta_m$, we have $u=[uu^T]\ones_m$, where $\ones_m$ is the all-ones vector in $\bR^m$. This observation allows to rewrite $F(u)$ as a homogeneous quadratic form:
\begin{equation}\label{exquadform}
F(u)=u^T\bar{Q}u,\,\,\bar{Q}=Q+[q\ones_m^T+\ones_mq^T].
\end{equation}
\par
Our goal is to construct an estimate $\widehat{g}(\zeta^{K})$ of $F(u)$, specifically, estimate of the form
\[
\widehat{g}(\zeta^{K})=\Tr (h\omega[\zeta^K])+\kappa
\]
where $\omega[\zeta^K]$ is the ``quadratic lifting'' of observation $\zeta^K$ (cf. \rf{liftedo1}):
\[
\omega[\zeta^K]={2\over K(K-1)}\sum_{1\leq j<j\leq M} \omega_{ij}[\zeta^K],\;\;\omega_{ij}[\zeta^{K}]=\half[\zeta_i\zeta_j^T+\zeta_j\zeta_i^T],\,1\leq i<j\leq K,
\]
and $h\in\bS^m$ and $\kappa\in\bR$ are the parameters of the estimate. To this end
\begin{itemize}
\item
we set $x(u)=uu^T$, with $X=\{uu^T:\;u\in U\}$, and specify a convex compact subset  $\X$ of the intersection of the ``symmetric matrix simplex'' $\bDelta^m\subset \bS^m$ (see \rf{mdeltad}) and the cone $\bS^m_+$ of positive semidefinite matrices such that
$X\subset \X\subset \E_X:=\bS^m$.
We put $\F=\E_F:=\bS^d$, and  $\M=\bDelta^d$,
thus $A\X A^T\subset \M\subset \E_M:=\bS^d$.
\item By Proposition \ref{discretelift}, $\F, \;\M\; $ and $\Phi(\cdot;\cdot)$, as defined in \rf{quaddphi}, form a regular data such that setting
    $
    M=\lfloor K/2\rfloor,
    $
for all $u\in U$  and  $h\in \bS^d$ it holds
\be
\begin{array}{c}
\ln\left(
\bE_{\zeta\sim P_{u}}\left\{\exp\{\langle h,\omega[\zeta^K]\rangle \}\right\}
\right)\leq \Phi_M\left(h;Auu^TA^T\right)\\
\left[\Phi_M(h;Z)=M\ln\left(\sum_{i,j}Z_{ij}\exp\{M^{-1}h_{ij}\}\right):\bS^d\times \bDelta^d\to\bR\right].\\
\end{array}
\ee{anatoli1d}
where $\langle h,w\rangle=\Tr(hw)$ is the Frobenius inner product on $\bS^d$.
\par
Observe that for $x\in \E_X$, $x\mapsto \A(x)=AxA^T$ is an affine mapping from $\X$ into $\M$, and setting
$$
G(x)=\langle \bar{Q},x\rangle:\E_X\to \bR,
$$
we get a linear functional on $\E_X$ such that
we ensure that
\[
G(uu^T)=\langle \bar{Q},uu^T\rangle=F(u).
\]
\end{itemize}
The relation $\Phi(0,z)=0\;\forall z\in \M$ being obvious,   Proposition \ref{discretelift} combines with Proposition \ref{corinf} to yield the following result.
\begin{proposition}\label{propquadform} In the situation in question, given $\epsilon\in(0,1)$, let $M=M(K)=\lfloor K/2\rfloor$, and let
\bse
\Psi_+(h,\alpha)&=&\max\limits_{x\in \X}\left[\alpha\Phi_M(h/\alpha,AxA^T)-\Tr(\bar{Q}x)\right]:\bS^d\times\{\alpha>0\}\to\bR,\\
\Psi_-(h,\alpha)&=&\max\limits_{x\in \X}\left[\alpha\Phi_M(-h/\alpha,AxA^T)+\Tr(\bar{Q}x)\right]:\bS^d\times\{\alpha>0\}\to\bR\\
\widehat{\Psi}_+(h)&:=&\inf_{\alpha>0}\left\{\Psi_+({h},{\alpha})+\alpha\ln(2/\epsilon)\right\}\\
&=&\max\limits_{x\in\X}\inf\limits_{\alpha>0}\left[\alpha\Phi_M(h/\alpha,AxA^T)-\Tr(\bar{Q}x)+\alpha\ln(2/\epsilon)\right]\\
&=&\max\limits_{x\in\X}\inf\limits_{\beta>0}\left[\beta\Phi_1(h/\beta,AxA^T)-\Tr(\bar{Q}x)+{\beta\over M}\ln(2/\epsilon)\right]\quad [\beta=M\alpha],\\
\widehat{\Psi}_-(h)&:=&\inf\limits_{\alpha>0}\left\{\Psi_-({h},{\alpha})+\alpha\ln(2/\epsilon)\right\}\\
&=&\max\limits_{x\in\X}\inf\limits_{\alpha>0}\left[\alpha\Phi_M(-h/\alpha,AxA^T)+\Tr(\bar{Q}x)+\alpha\ln(2/\epsilon)\right]\\
&=&\max\limits_{x\in\X}\inf\limits_{\beta>0}\left[\beta\Phi_1(-h/\beta,AxA^T)+\Tr(\bar{Q}x)+{\beta\over M}\ln(2/\epsilon)\right]\quad [\beta=M\alpha].
\ese
The functions  $\widehat{\Psi}_\pm$ are real valued and convex on $\bS^m$, and every candidate solution  $\bar{h}$ to the convex optimization problem
\be
\Opt=\min_h\left\{\widehat{\Psi}(h):=\half\left[\widehat{\Psi}_+(h)+\widehat{\Psi}_-(h)\right]\right\},
\ee{LQFwhynottoconsiderWW}
induces the estimate
$$
\widehat{g}_{\bar h}(\zeta^{K})=\Tr(\bar h\omega[\zeta^{K}])+\kappa(\bar h),\;\;\kappa(h)={\widehat{\Psi}_-(h) - \widehat{\Psi}_+(h)\over 2},
$$
 of the functional of interest {\rm(\ref{exquadform})} via observation $\zeta^{K}$
with $\epsilon$-risk on $U$ not exceeding $\bar \rho=\widehat{\Psi}(\bar h)$:
$$
\forall (u\in U): \Prob_{\zeta^K\sim P_u^K}\{|F(u)-\widehat{g}_{\bar h}(\zeta^K)|>\bar\rho\}\leq\epsilon.
$$
\end{proposition}
\subsubsection{Numerical illustration}
\begin{figure}
\begin{tabular}{cc}
\includegraphics[scale=0.55]{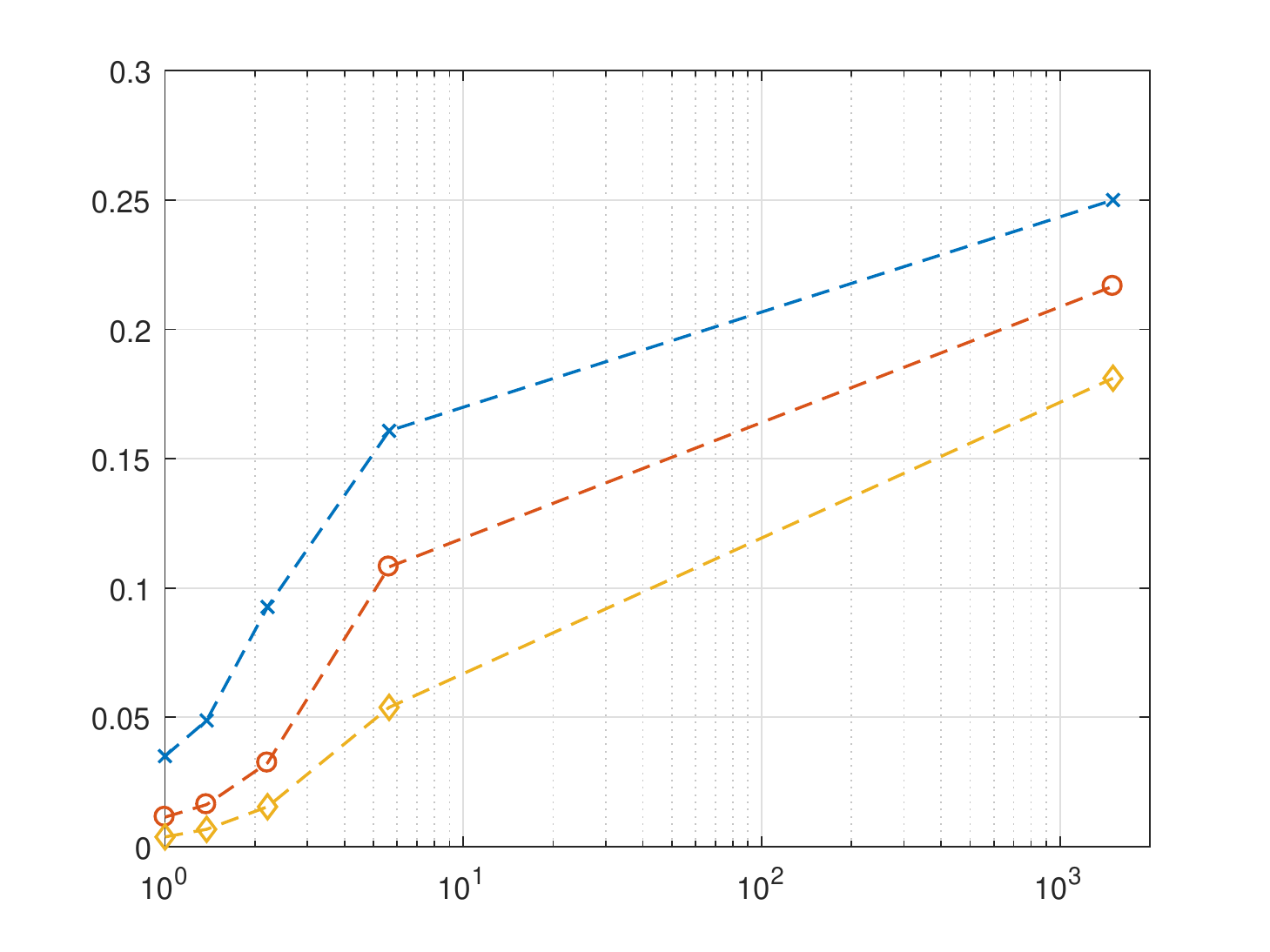}&\includegraphics[scale=0.55]{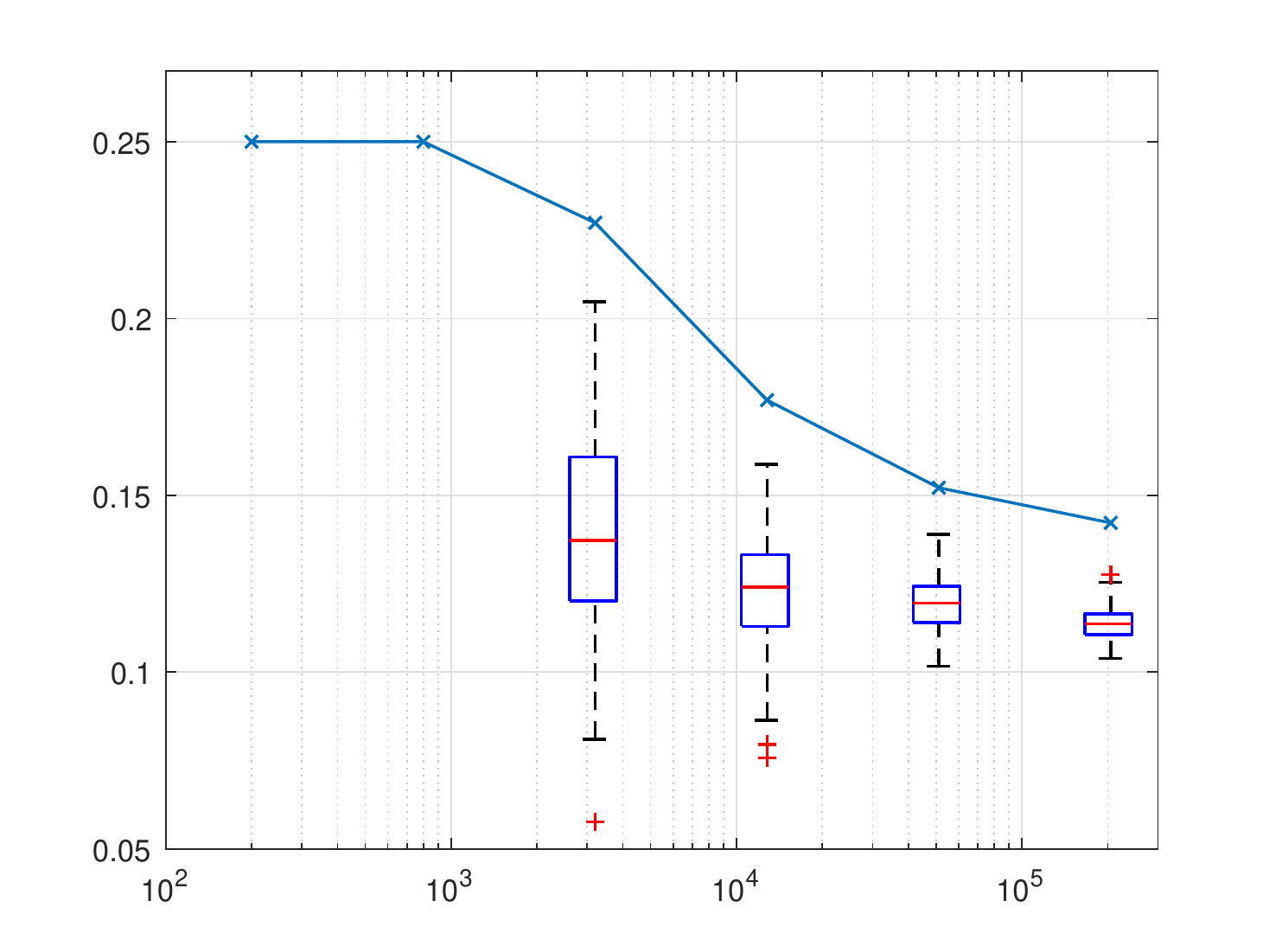}\\
$(a)$&$(b)$
\end{tabular}
\caption{\label{figdiscr}\small Estimation of ``independence defect.'' (a): Upper risk bound (value $\Opt$ in \rf{LQFwhynottoconsiderWW}) of linear estimate as a function of condition number $\mbox{\rm cond}(A)$; data for $K=2\cdot10^3,\;2\cdot10^4$ and $2\cdot10^5$. (b):  risk of linear estimation  as function of $K$ along with boxplots of empirical error distributions for $100$ simulations
($\theta=0.1$, $\hbox{\rm cond}(A)=39.2$).}
\end{figure}

To illustrate the above construction, consider the following problem: we observe independent across $k\leq K$ realizations $\zeta_k$ of discrete random variable $\zeta$ taking values $1,...,d$. The distribution $p\in\Delta^d$ of $\zeta$ is linearly parameterized by ``signal'' $u$ which itself is a probability distribution on ``discrete square'' $\Omega=\Xi\times\Xi$, $\Xi=\{1,...,m\}$:
$$
p_i=\sum\limits_{1\leq r,s\leq m} A_{p,rs}u_{rs},\,\,1\leq i\leq d.
$$
Here $A_{i,rs}\geq0$ are known coefficients such that $\sum_iA_{i,rs}=1$ for all $(r,s)\in\Omega$. Now, given two sets $I\subset \Xi$ and $J\subset\Xi$, consider the events  $\I=I\times\Xi\subset\Omega$ and $\cJ=\Xi\times J\subset\Omega$. Our objective is to quantify
the deviation of these events, the probability distribution on $\Omega$ being $u$, from independence, specifically, to estimate, via observations $\zeta_1,...,\zeta_K$, the quantity
$$
F_{IJ}(x)={\sum}_{(r,s)\in I\times J}u_{rs}-\left[{\sum}_{(r,s)\in I\times\Xi}u_{rs}\right]\left[{\sum}_{(r,s)\in \Xi\times J}u_{rs}\right]
$$
which is a quadratic function of $u$. In the experiments we report below, this estimation was carried out via a straightforward implementation of the construction presented earlier in this section. Our setup was as follows:
\begin{enumerate}
\item {We use $d=m^2$.} $d\times d$ column-stochastic ``sensing matrix'' $A$ \footnote{we identify the $m\times m$ ``discrete square'' $\Omega$ with $\{1,...,d\}$, which allows to treat a probability distribution $u$ on $\Omega$ as a vector from $\Delta_{d}$.} {corresponding to the ``mixed-noise observations'' \cite{lepski2016,lepski2017}} is generated according to $A=\theta I_{d}+(1-\theta)D$, with column-stochastic $d\times d$ matrix $D$, $\theta\in[0,1]$ being our control parameter. $D$ was selected at random, by normalizing columns of a $d\times d$ matrix with independent entries drawn from  the uniform distribution on $[0,1]$;
\item We set
\[
\X=\{x\in\bS^{d}:\;x_{rs,r's'}\geq0\,\forall r,s,r',s'\leq m,\, x\succeq0, \sum_{1\leq r,s,r',s'\leq m}x_{rs,r's'}=1\}
\]
which is the simplest convex outer approximation of the set $\{uu^T:u\in\Delta_{d}\}$.
\item We use $I=J=\{1,2,3\}\subset\Xi=\{1,2,...,8\}$, $\epsilon=0.01$,  $m=8$ (i.e., $d=64$).
\end{enumerate}
We present in Figure \ref{figdiscr} the results of experiments for $\theta$ taking values in $\{0.00,0.25,0.50,0,75,1.00\}$.
Other things being equal, the smaller $\theta$, the larger is the condition number $\hbox{\rm cond}(A)$ of the sensing matrix, and thus the larger is the (upper bound on the) risk of our estimate -- the optimal value of (\ref{LQFwhynottoconsiderWW}). Note that the variation of $F_{ij}$ over $X$ is exactly $1/2$, so the maximal risk is $\leq 1/4$. It {is worthy} to note that simple (if compared, e.g., to much more involved results of \cite{houdre2003exponential}) bounds in Proposition \ref{discretelift} for Laplace functional of order-2 $U$-statistics distribution
result in fairly good approximations of the risk of {our} estimate (cf. the boxplots of empirical distributions of the estimation error in the right plot of Figure \ref{figdiscr}).


\begin{thebibliography}{10}

\bibitem{bickel1988estimating}
P.~J. Bickel and Y.~Ritov.
\newblock Estimating integrated squared density derivatives: sharp best order
  of convergence estimates.
\newblock {\em Sankhy{\=a}: The Indian Journal of Statistics, Series A}, pages
  381--393, 1988.

\bibitem{Birge1981}
L.~Birg{\'e}.
\newblock Vitesses maximales de d{\'e}croissance des erreurs et tests optimaux
  associ{\'e}s.
\newblock {\em Zeitschrift f{\"u}r Wahrscheinlichkeitstheorie und verwandte
  Gebiete}, 55(3):261--273, 1981.

\bibitem{Birge1982}
L.~Birg\'e.
\newblock Sur un th\'eor\`eme de minimax et son application aux tests.
\newblock {\em Probab. Math. Stat.}, 3:259--282, 1982.

\bibitem{Birge1983}
L.~Birg{\'e}.
\newblock Approximation dans les espaces m{\'e}triques et th{\'e}orie de
  l'estimation.
\newblock {\em Zeitschrift f{\"u}r Wahrscheinlichkeitstheorie und verwandte
  Gebiete}, 65(2):181--237, 1983.

\bibitem{Birge2006}
L.~Birg{\'e}.
\newblock Model selection via testing: an alternative to (penalized) maximum
  likelihood estimators.
\newblock In {\em Annales de l'Institut Henri Poincare (B) Probability and
  Statistics}, volume~42, pages 273--325. Elsevier, 2006.

\bibitem{birge1995estimation}
L.~Birg{\'e} and P.~Massart.
\newblock Estimation of integral functionals of a density.
\newblock {\em The Annals of Statistics}, pages 11--29, 1995.

\bibitem{butucea2009adaptive}
C.~Butucea and F.~Comte.
\newblock Adaptive estimation of linear functionals in the convolution model
  and applications.
\newblock {\em Bernoulli}, 15(1):69--98, 2009.

\bibitem{butmez2011}
C.~Butucea and K.~Meziani.
\newblock Quadratic functional estimation in inverse problems.
\newblock {\em Statistical Methodology}, 8(1):31--41, 2011.

\bibitem{Caoetal2018}
Y.~Cao, A.~Nemirovski, Y.~Xie, V.~Guigues, and A.~Juditsky.
\newblock Change detection via affine and quadratic detectors.
\newblock {\em Electronic Journal of Statistics}, 12(1):1--57, 2018.

\bibitem{Don95}
D.~L. Donoho.
\newblock Statistical estimation and optimal recovery.
\newblock {\em The Annals of Statistics}, 22(1):238--270, 1994.

\bibitem{donoho1991geom2}
D.~L. Donoho and R.~C. Liu.
\newblock Geometrizing rates of convergence, ii.
\newblock {\em The Annals of Statistics}, pages 633--667, 1991.

\bibitem{donoho1991geom3}
D.~L. Donoho and R.~C. Liu.
\newblock Geometrizing rates of convergence, iii.
\newblock {\em The Annals of Statistics}, pages 668--701, 1991.

\bibitem{donoho1990minimax}
D.~L. Donoho, R.~C. Liu, and B.~MacGibbon.
\newblock Minimax risk over hyperrectangles, and implications.
\newblock {\em The Annals of Statistics}, pages 1416--1437, 1990.

\bibitem{donu1990}
D.~L. Donoho and M.~Nussbaum.
\newblock Minimax quadratic estimation of a quadratic functional.
\newblock {\em Journal of Complexity}, 6(3):290--323, 1990.

\bibitem{efromovich1996optimal}
S.~Efromovich and M.~Low.
\newblock On optimal adaptive estimation of a quadratic functional.
\newblock {\em The Annals of Statistics}, 24(3):1106--1125, 1996.

\bibitem{efromovich1994adaptive}
S.~Efromovich and M.~G. Low.
\newblock Adaptive estimates of linear functionals.
\newblock {\em Probability theory and related fields}, 98(2):261--275, 1994.

\bibitem{fan1991estimation}
J.~Fan.
\newblock On the estimation of quadratic functionals.
\newblock {\em The Annals of Statistics}, pages 1273--1294, 1991.

\bibitem{gayraud1999wavelet}
G.~Gayraud and K.~Tribouley.
\newblock Wavelet methods to estimate an integrated quadratic functional:
  Adaptivity and asymptotic law.
\newblock {\em Statistics \& probability letters}, 44(2):109--122, 1999.

\bibitem{GJN2015}
A.~Goldenshluger, A.~Juditsky, and A.~Nemirovski.
\newblock Hypothesis testing by convex optimization.
\newblock {\em Electronic Journal of Statistics}, 9(2):1645--1712, 2015.

\bibitem{cvx2014}
M.~Grant and S.~Boyd.
\newblock {\em The {\tt CVX} Users’ Guide. {Release} 2.1}, 2014.
\newblock \url{http://web.cvxr.com/cvx/doc/CVX.pdf}.

\bibitem{hasminskii1980some}
R.~Z. Hasminskii and I.~A. Ibragimov.
\newblock Some estimation problems for stochastic differential equations.
\newblock In {\em Stochastic Differential Systems Filtering and Control}, pages
  1--12. Springer, 1980.

\bibitem{houdre2003exponential}
C.~Houdr{\'e} and P.~Reynaud-Bouret.
\newblock Exponential inequalities, with constants, for u-statistics of order
  two.
\newblock In {\em Stochastic inequalities and applications}, pages 55--69.
  Springer, 2003.

\bibitem{huang1999nonparametric}
L.-S. Huang and J.~Fan.
\newblock Nonparametric estimation of quadratic regression functionals.
\newblock {\em Bernoulli}, 5(5):927--949, 1999.

\bibitem{ibragimov1988estimation}
I.~A. Ibragimov and R.~Z. Khas’~minskii.
\newblock Estimation of linear functionals in gaussian noise.
\newblock {\em Theory of Probability \& Its Applications}, 32(1):30--39, 1988.

\bibitem{ibragimov1987some}
I.~A. Ibragimov, A.~S. Nemirovskii, and R.~Khas’~minskii.
\newblock Some problems on nonparametric estimation in gaussian white noise.
\newblock {\em Theory of Probability \& Its Applications}, 31(3):391--406,
  1987.

\bibitem{JN2009}
A.~Juditsky and A.~Nemirovski.
\newblock Nonparametric estimation by convex programming.
\newblock {\em The Annals of Statistics}, 37(5a):2278--2300, 2009.

\bibitem{PartI}
A.~Juditsky and A.~Nemirovski.
\newblock Hypothesis testing via affine detectors.
\newblock {\em Electronic journal of statistics}, 10(2):2204--2242, 2016.

\bibitem{klemela2006sharp}
J.~Klemel{\"a}.
\newblock Sharp adaptive estimation of quadratic functionals.
\newblock {\em Probability theory and related fields}, 134(4):539--564, 2006.

\bibitem{klemela2001sharp}
J.~Klemela and A.~B. Tsybakov.
\newblock Sharp adaptive estimation of linear functionals.
\newblock {\em Annals of statistics}, pages 1567--1600, 2001.

\bibitem{laurent1997estimation}
B.~Laurent.
\newblock Estimation of integral functionals of a density and its derivatives.
\newblock {\em Bernoulli}, 3(2):181--211, 1997.

\bibitem{laurent2005adaptive}
B.~Laurent.
\newblock Adaptive estimation of a quadratic functional of a density by model
  selection.
\newblock {\em ESAIM: Probability and Statistics}, 9:1--18, 2005.

\bibitem{laurent2000adaptive}
B.~Laurent and P.~Massart.
\newblock Adaptive estimation of a quadratic functional by model selection.
\newblock {\em Annals of Statistics}, pages 1302--1338, 2000.

\bibitem{lepski2016}
O.~Lepski.
\newblock Some new ideas in nonparametric estimation.
\newblock {\em arXiv preprint arXiv:1603.03934}, 2016.

\bibitem{lepski2017}
O.~Lepski and T.~Willer.
\newblock Estimation in the convolution structure density model. part i: oracle
  inequalities.
\newblock {\em arXiv preprint arXiv:1704.04418}, 2017.

\bibitem{lepski1997optimal}
O.~V. Lepski and V.~G. Spokoiny.
\newblock Optimal pointwise adaptive methods in nonparametric estimation.
\newblock {\em The Annals of Statistics}, pages 2512--2546, 1997.

\bibitem{levit1975conditional}
B.~Y. Levit.
\newblock Conditional estimation of linear functionals.
\newblock {\em Problemy Peredachi Informatsii}, 11(4):39--54, 1975.

\end{thebibliography}

\appendix
\section{Proofs}
From now on, we use the notation
$$Z(u)=uu^T.
$$
\subsection{Proof of Proposition \ref{propGausslift0}}\label{sec:ppropGausslift}
Proposition \ref{propGausslift0} is nothing but \cite[Proposition 4.1.(i)]{Caoetal2018}; to make the paper self-contained, we reproduce the proof below. \par
We start with proving item (i) of Proposition.
\bigskip\par\noindent{\bf 1$^0$.} For any $\theta,h\in\bR^d,\;\Theta\in\bS^d_+$ and $H\in\bS^d$ such that  $-I\prec \Theta^{1/2}H\Theta^{1/2}\prec I$ we have
\be
\lefteqn{\Psi(h,H;\theta,\Theta):=\ln\left(\bE_{\zeta\sim\N(\theta,\Theta)}\left\{\exp\{h^T\zeta+\half \zeta^TH\zeta\}\right\}\right)}\nn
&=&\ln\left(\bE_{\xi\sim\N(0,I)}\left\{\exp\{h^T[\theta+\Theta^{1/2}\xi]+\half [\theta+\Theta^{1/2}\xi]^TH[\theta+\Theta^{1/2}\xi]\right\}\right) \nn
&=&-\half \ln\Det(I-\Theta^{1/2}H\Theta^{1/2})+h^T\theta+\half \theta^TH\theta+\half [H\theta+h]^T\Theta^{1/2}[I-\Theta^{1/2}H\Theta^{1/2}]^{-1}\Theta^{1/2}[H\theta+h]
\nn
&=&-\half \ln\Det(I-\Theta^{1/2}H\Theta^{1/2})
+\half [\theta;1]^T\hbox{\small $\left[\begin{array}{c|c}H&h\\\hline h^T&
\end{array}\right]$}[\theta;1]\nn
&&+\half [\theta;1]^T\left[[H,h]^T\Theta^{1/2}[I-\Theta^{1/2}H\Theta^{1/2}]^{-1}\Theta^{1/2}[H,h]\right][\theta;1].
\ee{56eq2}
 Observe that for $H\in \H_\gamma$ we have $\Theta^{1/2}[I-\Theta^{1/2}H\Theta^{1/2}]^{-1}\Theta^{1/2}=[\Theta^{-1}-H]^{-1}\preceq
[\Theta_*^{-1}-H]^{-1},$ so that \rf{56eq2} implies that for all $\theta\in\bR^d,\;\Theta\in\V,$ and $(h,H)\in\F$,
\be
\Psi(h,H;{\theta},\Theta)&\leq&-\half \ln\Det(I-\Theta^{1/2}H\Theta^{1/2})+\half
[\theta;1]^T\underbrace{\left[\hbox{\small $\left[\begin{array}{c|c}H&h\\\hline h^T&
\end{array}\right]$}+[H,h]^T[\Theta_*^{-1}-H]^{-1}[H,h]\right]}_{P[H,h]}[\theta;1]\nn
&=&-\half \ln\Det(I-\Theta^{1/2}H\Theta^{1/2})+\half \Tr(P[H,h]Z([\theta;1]))\nn
&=&-\half \ln\Det(I-\Theta^{1/2}H\Theta^{1/2})+\Gamma(h,H;Z([\theta;1])).
\ee{56eq22}
\bigskip\par\noindent{\bf 2$^0$.} We need the following
\begin{lemma}\label{lemlogdet} Let $\Theta_*$ be a $d\times d$ symmetric positive definite matrix, let $\delta\in[0,2]$, and let $\V$ be a closed convex subset of $\bS^d_+$ such that
\begin{equation}\label{suchthatweq1}
\Theta\in\V\Rightarrow\{\Theta\preceq\Theta_*\}\ \&\ \{\|\Theta^{1/2}\Theta_*^{-1/2}-I\|\leq\delta\}
\end{equation}
(cf. {\rm (\ref{56delta})}). Let also $\H^o:=\{H\in\bS^d:-\Theta_*^{-1}\prec H\prec \Theta_*^{-1}\}$. Then for all $(H,\Theta)\in\H^o\times\V$,
\begin{equation}\label{itholdsweq1}
-\half \ln\Det(I-\Theta^{1/2}H\Theta^{1/2})\leq \Upsilon(H;\Theta),
\end{equation}
where
\[
\Upsilon(H;\Theta)=-\half\ln\Det(I-\Theta_*^{1/2}H\Theta_*^{1/2}) +\half\Tr([\Theta-\Theta_*]H)+{\delta(2+\delta)\|\Theta_*^{1/2}H\Theta_*^{1/2}\|_F^2\over 2(1-\|\Theta_*^{1/2}H\Theta_*^{1/2}\|)}
\]
(here $\|\cdot\|$ is the spectral, and $\|\cdot\|_F$ - the Frobenius norm of a matrix).  \par
In addition, $\Upsilon(H,\Theta)$ is continuous function on $\H^o\times\V$ which is convex in $H\in H^o$ and concave (in fact, affine)
in $\Theta\in \V$
\end{lemma}
{\bf Proof.}   For $H\in\H^o$ and $\Theta\in\V$ fixed  we have
\[
\begin{array}{rcl}
\|\Theta^{1/2}H\Theta^{1/2}\|&=&\|[\Theta^{1/2}\Theta_*^{-1/2}][\Theta_*^{1/2}H\Theta_*^{1/2}][\Theta^{1/2}\Theta_*^{-1/2}]^T\|\\
&\leq& \|\Theta^{1/2}\Theta_*^{-1/2}\|^2\|\Theta_*^{1/2}H\Theta_*^{1/2}\|\leq \|\Theta_*^{1/2}H\Theta_*^{1/2}\|=:d(H)
\end{array}
\]
{with $d(H) <1$ for $H\in\H^o$}
(we have used the fact that $0\preceq\Theta\preceq\Theta_*$ implies $\|\Theta^{1/2}\Theta_*^{-1/2}\|\leq 1$).
Noting that $\|AB\|_F\leq\|A\|\|B\|_F$, {a similar computation} yields
\begin{equation}\label{56eq1235}
\|\Theta^{1/2}H\Theta^{1/2}\|_F\leq \|\Theta_*^{1/2}H\Theta_*^{1/2}\|_F=:D(H)
\end{equation}
Besides this, setting $F(X)=-\ln\Det(X):\inter\bS^d_+\to\bR$ and equipping $\bS^d$ with the Frobenius inner product, we have $\nabla F(X)=-X^{-1}$, so that with $R_0=\Theta_*^{1/2}H\Theta_*^{1/2}$, $R_1=\Theta^{1/2}H\Theta^{1/2}$, and $\Delta=R_1-R_0$, we have for properly selected $\lambda\in(0,1)$ and $R_\lambda=\lambda R_0+(1-\lambda)R_1$:
\bse
F(I-R_1)&=&F(I-R_0-\Delta)=F(I-R_0)+\langle \nabla F(I-R_\lambda),-\Delta\rangle =F(I-R_0)+\langle(I-R_\lambda)^{-1},\Delta\rangle\\
&=&F(I-R_0)+\langle I,\Delta\rangle +\langle (I-R_\lambda)^{-1}-I,\Delta\rangle.
\ese
We conclude that
\begin{equation}\label{56eqtru}
F(I-R_1)\leq F(I-R_0)+\Tr(\Delta)+\|I-(I-R_\lambda)^{-1}\|_F\|\Delta\|_F.
\end{equation}
Denoting by $\mu_i$ the eigenvalues of $R_\lambda$ and noting that $\|R_\lambda\|\leq \max[\|R_0\|,\|R_1\|]=d(H)$
we get $|\mu_i|\leq d(H)$. Therefore, the eigenvalues
\[\nu_i=1-{1\over 1-\mu_i}=-{\mu_i\over 1-\mu_i}
\] of
$I-(I-R_\lambda)^{-1}$ satisfy
\[|\nu_i|\leq |\mu_i|/(1-\mu_i)\leq |\mu_i|/(1-d(H)),
\] whence $$\|I-(I-R_\lambda)^{-1}\|_F\leq \|R_\lambda\|_F/(1-d(H)).$$
Noting that $\|R_\lambda\|_F\leq\max[\|R_0\|_F,\|R_1\|_F]\leq D(H)$, see (\ref{56eq1235}), we conclude that \[\|I-(I-R_\lambda)^{-1}\|_F\leq D(H)/(1-d(H)),\] and when substituting into (\ref{56eqtru}) we get
\begin{equation}\label{56eqtru1}
F(I-R_1)\leq F(I-R_0)+\Tr(\Delta)+D(H)\|\Delta\|_F/(1-d(H)).
\end{equation}
Furthermore, because by (\ref{56delta}) the matrix $D =\Theta^{1/2}\Theta_*^{-1/2}-I$ satisfies $\|D\|\leq\delta$,
$$
\Delta=\underbrace{\Theta^{1/2}H\Theta^{1/2}}_{R_1}-\underbrace{\Theta_*^{1/2}H\Theta_*^{1/2}}_{R_0}=
(I+D)R_0(I+D^T)-R_0=DR_0+R_0D^T+DR_0D^T.
$$
Consequently,
\bse
\|\Delta\|_F&\leq& \|DR_0\|_F+\|R_0D^T\|_F+\|DR_0D^T\|_F\leq [2\|D\|+\|D\|^2]\|R_0\|_F\\&\leq& \delta(2+\delta)\|R_0\|_F=
\delta(2+\delta)D(H).
\ese
This combines with (\ref{56eqtru1}) and the relation
\[
\Tr(\Delta)=\Tr(\Theta^{1/2}H\Theta^{1/2}-\Theta_*^{1/2}H\Theta_*^{1/2})=\Tr([\Theta-\Theta_*]H)
\] to yield
$$
F(I-R_1)\leq F(I-R_0)+\Tr([\Theta-\Theta_*]H)+{\delta(2+\delta)\over1-d(H)}\|\Theta_*^{1/2}H\Theta_*^{1/2}\|_F^2,
$$
and we arrive at (\ref{itholdsweq1}). It remains to prove that $\Upsilon(H;\Theta)$ is convex-concave and continuous on $\H^o\times\V$. The only component of this
claim which is not completely evident is convexity of the function in $H\in\H^o$. To see that it is indeed the case, note that $\ln\Det(\cdot)$ is concave on the interior of the semidefinite cone, function $f(u,v)={u^2\over 1-v}$ is convex and nondecreasing in $u,v$ in the convex domain $\Pi=\{(u,v):u\geq0,v<1\}$, and the function ${\|\Theta_*^{1/2}H\Theta_*^{1/2}\|_F^2\over 1-\|\Theta_*^{1/2}H\Theta_*^{1/2}\|}$ is obtained from $f$ by convex substitution of variables $H\mapsto(\|\Theta_*^{1/2}H\Theta_*^{1/2}\|_F,\|\Theta_*^{1/2}H\Theta_*^{1/2}\|)$ mapping
$\H^o$ into $\Pi$.
\qed
\bigskip\par\noindent{\bf 3$^0$.} Combining (\ref{itholdsweq1}), (\ref{56eq22}), \rf{PhisubG} and the origin of $\Psi$, see (\ref{56eq2}), we conclude
 that for all $(\theta,\Theta)\in \bR^d\times\V$ and $(h,H)\in\F=\bR^d\times \H_\gamma)$,
$$
\ln\left(\bE_{\zeta\sim\N(\theta,\Theta)}\left\{\exp\{h^T\zeta+\half \zeta^TH\zeta\}\right\}\right)\leq\Phi(h,H;\Theta,Z([\theta;1])).
$$
To complete the proof of (i), all we need is to verify  the claim that $\F,\M^+,\Phi$ is regular data, which boils down to checking that $\Phi:\F\times\M^+\to\bR$ is {continuous and} convex-concave. Let us verify convexity-concavity and continuity. Recalling that $\Upsilon(H;\Theta):
\F\times\V\to\bR$ indeed is convex-concave and continuous, the verification in question reduces to checking that $\Gamma(h,H;Z)$ is convex-concave and continuous on $(\bR^d\times \H_\gamma)\times \Z^+$. Continuity and concavity in $Z$ being evident, all we need to prove is that whenever $Z\in\Z^+$, the function $\Gamma_Z(h,H):=\Gamma(h,H;Z)$ is convex in
$(h,H)\in\F=\bR^d\times \H_\gamma$. By the Schur Complement Lemma, we have
$$
\G:=\{(h,H,G):\; G\succeq P[H,h]\}=\left\{(h,H,G):\left[\begin{array}{c|c}G-\hbox{\small $\left[\begin{array}{c|c}H&h\\\hline h^T&
\end{array}\right]$}&[H,h]^T\cr\hline
[H,h]&\Theta_*^{-1}-H\cr\end{array}\right]\succeq0\right\},
$$
implying that $\G$ is convex. Now, since $Z\succeq0$ due to $Z\in \Z^+\subset\bS^{m+1}_+$, we have
$$
\{(h,H,\tau):\;(h,H)\in \H_\gamma,\;\tau\geq\Gamma_Z(h,H)\}=\{(h,H,\tau):(h,H)\in \H_\gamma,\,
\exists G: G\succeq {P}[H,h],\;2\tau\geq \Tr(ZG)\},
$$
and because $\G$ is convex, so is the epigraph of $\Gamma_Z$, as claimed. Item (i) of Proposition \ref{propGausslift0} is proved.
\bigskip\par\noindent
{\bf 4$^0$.} It remains to verify  item (ii) of Proposition \ref{propGausslift0} stating that $\Phi$ is coercive in $h,H$. Let $\Theta\in\V$, $Z\in\Z^+$, and $(h_i,H_i)\in\bR^d\times \H_\gamma$ with $\|(h_i,H_i)\|\to\infty$ as $i\to\infty$,
and let us prove that $\Phi(h_i,H_i;\Theta,Z)\to\infty$. Looking at the expression for $\Phi(h_i,H_i;\Theta,Z)$, it is immediately seen that all terms in this expression, except for the terms coming from $\Gamma(h_i,H_i;Z)$,
remain bounded as $i$ grows, so that all we need to verify  is that $\Gamma(h_i,H_i;Z)\to \infty$ as $i\to\infty$. Observe that the sequence $\{H_i\}_i$ is bounded due to $H_i\in \H_\gamma$, implying that $\|h_i\|_2\to\infty$ as $i\to\infty$. Denoting by $e$ the last basic orth of $\bR^{d+1}$ and taking into account that the matrices $[\Theta_*^{-1}-H_i]^{-1}$ satisfy $\alpha I_d\preceq [\Theta_*^{-1}-H_i]^{-1}\preceq \beta I_d$  for some positive $\alpha,\beta$ due to $H_i\in\H_\gamma$, observe that
$$
\begin{array}{l}
\underbrace{\left[\left[\begin{array}{c|c}H_i&h_i\cr\hline h_i^T&\end{array}\right]+
\left[H_i,h_i\right]^T[\Theta_*^{-1}-H_i]^{-1}\left[H_i,h_i\right]\right]}_{P_i}
=
\underbrace{\left[h_i^T[\Theta_*^{-1}-H_i]^{-1}h_i\right]}_{\alpha_i\|h_i\|_2^2}ee^T+R_i,\\
\end{array}
$$
where $\alpha_i\geq\alpha>0$ and $\|R_i\|_F\leq C(1+\|h_i\|_2)$. As a result,
\bse
\Gamma(h_i,H_i;Z)&\geq&\Tr(ZP_i)=\Tr(Z[\alpha_i\|h_i\|_2^2ee^T+R_i])\\
&\geq&\alpha_i\|h_i\|_2^2
\underbrace{\Tr(Zee^T)}_{=Z_{m+1,m+1}=1}-\|Z\|_F\|R_i\|_F
\geq
\alpha\|h_i\|_2^2-C(1+\|h_i\|_2)\|Z\|_F,\\
\ese
and the concluding quantity tends to $\infty$ as $i\to\infty$ due to $\|h_i\|_2\to\infty$, $i\to\infty$. \qed
\subsection{Proof of Proposition \ref{discretelift}}\label{sec:pprop2}
Continuity and convexity-concavity of $\Phi$ and $\Phi_M$ are obvious. Let us verify relations \rf{excovering}. Let us fix $\mu\in \Delta^d$, and let
and $\zeta^K\sim P^K_\mu$.
Let us denote $\S_K$ the set of all permutations $\sigma$ of $\{1,...,K\}$, and let
$$
\omega^\sigma[\zeta^{K}]={1\over M}\sum_{k=1}^M\omega_{\sigma_{2k-1}\sigma_{2k}}[\zeta^{K}],\;\;\sigma\in\S_K.
$$
By the symmetry argument we clearly have
$$
\sum_{\sigma\in\S_K}\sum_{k=1}^M\omega_{\sigma_{2k-1}\sigma_{2k}}[\zeta^K]=N\sum_{1\leq i\neq j\leq K}\zeta_{i}\zeta_{j}^T=
2N\sum_{1\leq i<j\leq K}\omega_{ij}[\zeta^K],
$$
where $N$ is the number of permutations $\sigma\in\S_K$ such that a particular pair $(i,j)$, $1\leq i\neq j\leq K$ is met among the pairs $(\sigma_{2k-1},\sigma_{2k})$, $1\leq k\leq M$. Comparing the total number of $\omega_{ij}$-terms in the left and the right hand sides of the latter equality, we get $\Card(\S_K)M=NK(K-1)$, which combines with the equality itself to imply that
\be
{2\over K(K-1)}\sum_{1\leq i<j\leq K}\omega_{ij}[\zeta^K]={1\over \sCard(\S_K)}\sum_{\sigma\in\S_K}{1\over M}\sum_{k=1}^M\omega_{\sigma_{2k-1}\sigma_{2k}}[\zeta^K]={1\over \sCard(\S_K)}\sum_{\sigma\in \S_K}\omega^\sigma[\zeta^{K}].
\ee{sigma_2k}
 Let $\sigma_{\tid}$ be the identity permutation of $1,...,K$. Due to \rf{sigma_2k} we have
\be
\bE_{\zeta^{K}\sim P_\mu^K}\left\{\exp\{\Tr(H\omega[\zeta^{K}])\}\right\}&=&
\bE_{\zeta^{K}\sim P_\mu^K}\left\{\exp\left\{{1\over \sCard(\S_K)}\sum_{\sigma\in\S_K}\Tr(H\omega^\sigma[\zeta^{K}])\right\}\right\}\nn
\hbox{\ [by the H\"{o}lder inequality]}&\leq& \prod\limits_{\sigma\in\S_K} \left[\bE_{\zeta^{K}\sim P_\mu^K}\left\{\exp\left\{\Tr(H\omega^\sigma[\zeta^{K}])\right\}\right\}\right]^{1/\sCard(\S_K)}
\nn
\hbox{[because $\omega^\sigma[\zeta^{K}]$ are equally distributed $\forall \sigma$]}&=&
\bE_{\zeta^{K}\sim P_\mu^K} \left\{\exp\{\Tr(H\omega^{\sigma_{\tid}}[\zeta^{K}])\}\right\}\nn
\hbox{\ [by definition of $\omega^{\sigma_{\tid}}[\cdot]$]}&=&\bE_{\zeta^{K}\sim P_\mu^K} \left\{\prod_{k=1}^M\exp\left\{{1\over M}\Tr(H\omega_{2k-1,2k}[\zeta^K])\right\}\right\}\nn
\hbox{\ [since $\zeta_1,...,\zeta_{K}$ are i.i.d.]}&=&\left[\bE_{\zeta^{2}\sim P_\mu^2}\left\{\exp\{\Tr((H/M)\omega_{12}[\zeta^2])\}\right\}\right]^M.
\ee{combi}
The distribution of the random variable $\omega_{12}[\zeta^2]=\half[\zeta_1\zeta_2^T+\zeta_2\zeta_1^T]$, $\zeta^2\sim P_\mu^2$, clearly is $P_{Z[\mu]}$, so that
\bse
\ln\left(\bE_{\zeta^{2}\sim P_\mu^2}\left\{\exp\{\Tr((H/M)\omega_{12}[\zeta^2])\}\right\}\right)&=&
\ln\left(\bE_{W\sim P_{Z[\mu]}}\left\{\exp\{\Tr((H/M)W)\}\right\}\right)\\&=&\ln\left(\sum_{i,j=1}^d e^{M^{-1}H_{ij}}\mu_i\mu_j\right)=\Phi(H/M;Z(\mu)).
\ese
The latter relation combines with \rf{combi} to imply \rf{excovering}. \qed

\subsection{Proof of Proposition \ref{corinf}}\label{proof:corinf}

Let us first verify the identities \rf{psisareequala} and \rf{psisareequalb}. The function
$$
\Theta(f,\alpha;x)=\alpha\Phi(f/\alpha,\A(x))-G(x)+\alpha\ln(2/\epsilon):\,\F^+\times \X\to\bR
$$
is convex-concave and continuous, and $\X$ is compact. Hence, by Sion-Kakutani Theorem,
$$
\begin{array}{rcl}
\widehat{\Psi}_+(f)&:=&\inf_{\alpha}\left\{\Psi_+({f},{\alpha})+\alpha\ln(2/\epsilon):\alpha>0,\;(f,\alpha)\in\F^+\right\}\\
&=&\inf_{\alpha>0,(f,\alpha)\in\F^+}\max_{x\in\X} \Theta(f,\alpha;x)=\sup_{x\in\X}\inf_{\alpha>0,(f,\alpha)\in\F^+} \Theta(f,\alpha;x)\\
&=&\sup_{x\in\X}\inf_{\alpha>0,(f,\alpha)\in\F^+}\left[\alpha\Phi(f/\alpha,\A(x))-G(x)+\alpha\ln(2/\epsilon)\right],\\
\end{array}
$$
as required in \rf{psisareequala}. As we know,  $\Psi_+(f,\alpha)$ is a real-valued continuous function on $\F^+$, so that $\widehat{\Psi}_+$ is convex on $\E_F$, provided that the function is real-valued. Now, let $\bar{x}\in\X$, and let $\psi$ be a subgradient of $\phi(f)=\Phi(f;\A(\bar{x}))$ taken at $f=0$. For $f\in\E_F$ and all $\alpha>0$ such that $(f,\alpha)\in\F^+$ we have
$$
\begin{array}{rcl}
\Psi_+(f,\alpha)&\geq& \alpha\Phi(f/\alpha;\A(\bar{x}))-G(\bar{x})+\alpha\ln(2/\epsilon)\\
&\geq& \alpha[\Phi(0;\A(\bar{x}))+\langle \psi, f/\alpha\rangle]-G(\bar{x})+\alpha\ln(2/\epsilon)\geq \langle \psi, f\rangle-G(\bar{x})\\
\end{array}
$$
(we have used (\ref{Phisatisfies})). Therefore, $\Psi_+(f,\alpha)$ is below bounded on the set $\{\alpha>0:f/\alpha\in\F\}$. In addition, this set is nonempty, since $\F$ contains a neighbourhood of the origin. Thus, $\widehat{\Psi}_+$ is real-valued and convex on $\E_F$.
 Verification of \rf{psisareequalb} and of the fact that $\widehat{\Psi}_-(f)$ is a real-valued convex function on $\E_F$ is completely similar.
\par Now, given a feasible solution $(\bar{f},\bar{\varkappa},\bar{\rho})$ to (\ref{sssuchthatcor}), let us select somehow $\widetilde{\rho}>\bar{\rho}$.
Taking into account the definition of $\widehat{\Psi}_\pm$, we can find $\bar{\alpha}$ and $\bar{\beta}$ such that
$$
\begin{array}{l}
(\bar{f},\bar{\alpha})\in\F^+\;\mbox{and}\; \Psi_+(\bar{f},\bar{\alpha})+\bar{\alpha}\ln(2/\epsilon)\leq\widetilde{\rho}-\bar{\varkappa},\\
(\bar{f},\bar{\beta})\in\F^+\;\mbox{and}\;\Psi_-(\bar{f},\bar{\beta})+\bar{\beta}\ln(2/\epsilon)\leq\widetilde{\rho}+\bar{\varkappa},\\
\end{array}
$$
implying that the collection $(\bar{f},\bar{\alpha},\bar{\beta},\bar{\varkappa},\widetilde{\rho})$ is a feasible solution to (\ref{sssuchthat}).
We need the following statement.
\begin{lemma}\label{lemlemlem} Given $\epsilon\in(0,1)$, let $\bar{f}$, $\bar{\alpha}$, $\bar{\beta}$, $\bar{\varkappa}$, $\widetilde{\rho}$
be a feasible solution to the system of convex constraints
\be
\begin{array}{l}
(a)~~~~~~({f},{\alpha})\in\F^+,\;{\alpha}\ln(\epsilon/2)\geq\Psi_+({f},{\alpha})-{\rho}+{\varkappa},\\
(b)~~~~~~({f},{\beta})\in\F^+,\;{\beta}\ln(\epsilon/2)\geq\Psi_-({f},{\beta})-{\rho}-{\varkappa},
\end{array}
\ee{sssuchthat}
in variables $f$, $\alpha$, $\beta$, $\rho$, $\varkappa$.
Then the $\epsilon$-risk of the estimate
$
\widehat{g}(\omega)=\langle \bar{f},\omega\rangle+\bar{\varkappa},
$
is at most $\widetilde{\rho}$.
\end{lemma}
\paragraph{Proof.}
Let $\epsilon\in(0,1)$, $\bar{f}$, $\bar{\alpha}$, $\bar{\beta}$, $\bar{\varkappa}$, $\widetilde{\rho}$  satisfy the premise of Lemma, and let $x\in X,P$ satisfy (\ref{123cond21}). We have
\bse
\Prob_{\omega\sim P}\{\widehat{g}(\omega)>G(x)+\widetilde{\rho}\}
&\leq &\left[\int{\rm e}^{\langle\bar{f},\omega\rangle/ \bar{\alpha}}P(d\omega)\right]
{\rm e}^{-{G(x)+\widetilde{\rho}-\bar{\varkappa}\over\bar{\alpha}}}
\leq {\rm e}^{\Phi\left({\bar{f}/ \bar{\alpha}},\A(x)\right)}{\rm e}^{-{G(x)+\widetilde{\rho}-\bar{\varkappa}\over \bar{\alpha}} }.
\ese
Thus
\bse
\bar{\alpha}\ln\left(\Prob_{\omega\sim P}\{\widehat{g}(\omega)>G(x) +\widetilde{\rho}\}\right)&\leq& \bar{\alpha}\Phi(\bar{f}/\bar{\alpha},\A(x))-G(x) -\widetilde{\rho}+\bar{\varkappa}\\
\hbox{\ [by definition of $\Psi_+$ and due to $x\in X$]}&\leq& \Psi_+(\bar{f},\bar{\alpha})-\widetilde{\rho}+\bar{\varkappa}\\
\hbox{\ [by (\ref{sssuchthat}.a)]}&\leq&\bar{\alpha}\ln(\epsilon/2),
\ese
and we conclude that
\[
\Prob_{\omega\sim P}\{\widehat{g}(\omega)>G(x)+\widetilde{\rho}\}\leq \epsilon/2.
\]
Similarly,
\bse
\Prob_{\omega\sim P}\{\widehat{g}(\omega)<G(x)-\widetilde{\rho}\}
&\leq& \left[\int{\rm e}^{-{\langle \bar{f},\omega\rangle/\bar{\beta}}}P(d\omega)\right]{\rm e}^{-{-G(x) +\widetilde{\rho}+\bar{\varkappa}\over \bar{\beta}}}
\leq {\rm e}^{\Phi\left(-{\bar{f}/ \bar{\beta}},\A(x)\right)}{\rm e}^{{G(x)-\widetilde{\rho}-\bar{\varkappa}\over\bar{\beta}}},\\
\ese
whence
\bse
\bar{\beta}\ln\left(\Prob_{\omega\sim P}\{\widehat{g}(\omega)<G(x) -\widetilde{\rho}\}\right)&\leq& \bar{\beta}\Phi(-\bar{f}/\bar{\beta},\A(x))+G(x) -\widetilde{\rho}-\bar{\varkappa}\\
\hbox{\ [by definition of $\Psi_-$ and due to $x\in X$]}&\leq& \Psi_-(\bar{f},\bar{\beta})-\widetilde{\rho}-\bar{\varkappa}\\
\hbox{\ [by (\ref{sssuchthat}.b)]}&\leq&\bar{\beta}\ln(\epsilon/2),
\ese
so that\[
\Prob_{\omega\sim P}\{\widehat{g}(\omega)<G(x) -\widetilde{\rho}\}\leq \epsilon/2.
\eqno{\hbox{\qed}}
\]
When invoking Lemma \ref{lemlemlem},  we get
$$
\Prob_{\omega\sim P} \left\{\omega: |\widehat{g}(\omega)-G(x)|>\widetilde{\rho}\right\}\leq\epsilon
$$
for all $(x\in X,P\in\P)$ satisfying (\ref{123cond21}). Since $\widetilde{\rho}$ can be selected arbitrarily close to $\bar{\rho}$, $\widehat{g}(\cdot)$ indeed is a $(\bar{\rho},\epsilon)$-accurate estimate. \qed

\subsection{Proof of Proposition \ref{melmelmel}}\label{proof:melmelmel}
Under the premise of the proposition, let us fix ${u}\in U$, $v\in V$, so that ${x}:=(v,Z(u):=[u;1][u;1]^T)\in X$. Denoting by $P=P_{u,v}$ the distribution of $\omega:=(\zeta,\zeta\zeta^T)$ with $\zeta\sim\N(A[{u};1],{M(v)})$, and
invoking \rf{itholds123}, we see that for just defined $(x,P)$, relation  (\ref{123cond21}) takes place. Applying Proposition \ref{corinfrepeated}, we conclude that
$$
\Prob_{{\zeta^K\sim[\N(A[u;1],M(v))]^K}}\left\{|\widehat{g}(\zeta^K)-G(x)|>\bar{\rho}\right\}\leq\epsilon.
$$
It remains to note that by construction it holds
$$
G\left(x=(v,Z([u;1]))\right)=q^Tv+\Tr(QZ([u;1]))=q^Tv+\Tr(Q[u;1][u;1]^T)=q^Tv+[u;1]^TQ[u,1]=F(u,v).
$$
The ``in particular'' part of proposition is immediate -- with $\rho$ and $\varkappa$ given by (\ref{varkapparho}),
$h,H,\rho,\varkappa$ clearly satisfy (\ref{sssuchthatthat}). \qed
\subsection{Proof of Proposition \ref{propconsistency}}\label{NEwProof}
\paragraph{1$^0$.} By {A.2}, the columns of $(d+1)\times (m+1)$ matrix $B$, see \rf{bmat}, are linearly independent, so that we can
find $(m+1)\times (d+1)$ matrix $C$ such that $CB=I_{m+1}$.
Let us define $(\bar{h},\bar{H})\in\bR^d\times \bS^d$ from the relation
\begin{equation}\label{eqbarHplus}
\barHplus=2[C^TQC]^o,
\end{equation}
where for $(d+1)\times(d+1)$ matrix $S$, $S^o$ is the matrix obtained from $S$ by replacing the entry $S_{d+1,d+1}$ with zero.
\paragraph{2$^0$.} Let us fix $\epsilon\in(0,1)$. Setting
$$
\rho_K=\half\left[\widehat{\Psi}^{K}_+(\bar{h},\bar{H})+\widehat{\Psi}^{K}_-(\bar{h},\bar{H})\right]\\
$$
and invoking Proposition \ref{melmelmel}, all we need to prove is that in the case of {A.1-2} one has
\begin{equation}\label{2onehas}
\lim\sup_{K\to\infty}\left[\widehat{\Psi}_+^{K}(\bar{h},\bar{H})+\widehat{\Psi}_-^{K}(\bar{h},\bar{H})\right]\leq0.
\end{equation}
To this end note that in our current situation, \rf{PhisubG} and (\ref{newnewnewnewnew}) simplify to
{\small$$
\begin{array}{l}
\Phi(h,H;Z)=-\half\ln\Det(I-\Theta_*^{1/2}H\Theta_*^{1/2})+\half\Tr\bigg(Z\underbrace{\left(B^T\left[\hbox{\small$\left[\begin{array}{c|c}H&h\cr\hline h^T&\end{array}\right]+
\left[H,h\right]^T[\Theta_*^{-1}-H]^{-1}\left[H,h\right]$}\right]B\right)}_{P[H,h]}\bigg),\\
\widehat{\Psi}_+^{K}(h,H)
=\inf\limits_{\alpha}\left\{\max\limits_{Z\in\Z}\left[\alpha\Phi(h/\alpha,H/\alpha;Z)-\Tr(QZ)+K^{-1}\alpha\ln(2/\epsilon)\right]:
\alpha>0,
-\gamma\alpha\Theta_*^{-1}\preceq H\preceq \gamma\alpha\Theta_*^{-1}\right\},\\
\widehat{\Psi}_-^{K}(h,H)
=\inf\limits_{\alpha}\left\{\max\limits_{Z\in\Z}\left[\alpha\Phi(-h/\alpha,-H/\alpha;Z)+\Tr(QZ)+K^{-1}\alpha\ln(2/\epsilon)\right]:
\alpha>0,
-\gamma\alpha\Theta_*^{-1}\preceq H\preceq \gamma\alpha\Theta_*^{-1}\right\}.\\
\end{array}
$$}
\noindent
Hence,
{\small
\begin{equation}\label{retro}
\begin{array}{ll}
\multicolumn{2}{l}{\left[\widehat{\Psi}_+^K(\bar{h},\bar{H})+\widehat{\Psi}_-^K(\bar{h},\bar{H})\right]
\leq \inf\limits_{\alpha}\bigg\{\max\limits_{Z_1,Z_2\in\Z}\bigg[\alpha\Phi(\bar{h}/\alpha,\bar{H}/\alpha;Z_1)
-\Tr(QZ_1)+\Phi(-\bar{h}/\alpha,-\bar{H}/\alpha;Z_2)+\Tr(QZ_2)}\\
&\multicolumn{1}{r}{+2K^{-1}\alpha\ln(2/\epsilon)\bigg]:\alpha>0,-\gamma\alpha\Theta_*^{-1}\preceq \bar{H}\preceq \gamma\alpha\Theta_*^{-1}\bigg\}}
\\
&=\inf\limits_{\alpha}\max\limits_{Z_1,Z_2\in\Z}\bigg\{-\half\alpha\ln\Det\left(I-[\Theta_*^{1/2}\bar{H}\Theta_*^{1/2}]^2/\alpha^2\right)+2K^{-1}\alpha\ln(2/\epsilon)+\Tr(Q[Z_2-Z_1])\\
&\multicolumn{1}{r}{+\half\left[
\alpha\Tr\left(Z_1P[\bar{H}/\alpha,\bar{h}/\alpha]\right)+\alpha\Tr\left(Z_2P[-\bar{H}/\alpha,-\bar{h}/\alpha]\right)\right]:\alpha>0,-\gamma\alpha\Theta_*^{-1}\preceq \bar{H}\preceq \gamma\alpha\Theta_*^{-1}\bigg\}}\\
&=\inf\limits_{\alpha}\max\limits_{Z_1,Z_2\in\Z}\bigg\{-\half\alpha\ln\Det\left(I-[\Theta_*^{1/2}\bar{H}\Theta_*^{1/2}]^2/\alpha^2\right)+2K^{-1}\alpha\ln(2/\epsilon)\\
&\multicolumn{1}{r}{+\underbrace{\Tr(Q[Z_2-Z_1])+\half\Tr([Z_1-Z_2]B^T\barHplus B)}_{T[Z_1,Z_2]}
+\half\Tr\left(Z_1B^T[\bar{H},\bar{h}]^T[\alpha\Theta_*^{-1}-\bar{H}]^{-1}[\bar{H},\bar{h}]B\right)}\\
&\multicolumn{1}{r}{+\half\Tr\left(Z_2B^T[\bar{H},\bar{h}]^T
[\alpha\Theta_*^{-1}+\bar{H}]^{-1}[\bar{H},\bar{h}]B\right):\alpha>0,-\gamma\alpha\Theta_*^{-1}\preceq \bar{H}\preceq \gamma\alpha\Theta_*^{-1}\bigg\}}
\end{array}
\end{equation}}
\noindent By (\ref{eqbarHplus}) we have
$
\half B^T\barHplus B=B^T[C^TQC+J]B,
$ where the only nonzero entry, if any, in the $(d+1)\times(d+1)$ matrix $J$ is {$J_{d+1,d+1}$}. Due to the structure of $B$, see  \rf{bmat}, we conclude that the only nonzero element, if any, in $\bar{J}=B^TJB$ is {$\bar{J}_{m+1,m+1}$}, and that
$$
\half B^T\barHplus B=(CB)^TQ(CB)+\bar{J}=Q+\bar{J}
$$
(recall that $CB=I_{m+1}$). Now, whenever $Z\in \Z$, one has $Z_{m+1,m+1}=1$, whence
$$
\half\Tr([Z_1-Z_2] B^T\barHplus B)=\Tr([Z_1-Z_2]Q)+\Tr([Z_1-Z_2]\bar{J})=\Tr([Z_1-Z_2]Q),
$$
implying that the quantity $T[Z_1,Z_2]$ in (\ref{retro}) is zero, provided $Z_1,Z_2\in \Z$. Consequently, (\ref{retro}) becomes
{\small\begin{equation}\label{retro1}
\begin{array}{ll}
\left[\widehat{\Psi}_+^K(\bar{h},\bar{H})+\widehat{\Psi}_-^K(\bar{h},\bar{H})\right]
&\leq\inf\limits_{\alpha}\max\limits_{Z_1,Z_2\in\Z}\bigg\{-\half\alpha\ln\Det\left(I-[\Theta_*^{1/2}\bar{H}\Theta_*^{1/2}]^2/\alpha^2\right)+2K^{-1}\alpha\ln(2/\epsilon)\\
&\multicolumn{1}{r}{
+\half\Tr\left(Z_1B^T[\bar{H},h][\alpha\Theta_*^{-1}-\bar{H}]^{-1}[\bar{H},\bar{h}]^TB\right)}\\
&\multicolumn{1}{r}{+\half\Tr\left(Z_2B^T[\bar{H},\bar{h}]^T
[\alpha\Theta_*^{-1}+\bar{H}]^{-1}[\bar{H},\bar{h}]B\right):\alpha>0,-\gamma\alpha\Theta_*^{-1}\preceq \bar{H}\preceq \gamma\alpha\Theta_*^{-1}\bigg\}}
\\
\end{array}
\end{equation}}
Now, for appropriately selected  independent of $K$ real $c>0$ we have for $\alpha\geq c$:
\[
-\half\alpha\ln\Det\left(I-[\Theta_*^{1/2}\bar{H}\Theta_*^{1/2}]^2/\alpha^2\right)\leq c/\alpha,
\]
and
\[
\half\Tr\left(Z_1B^T[\bar{H},\bar{h}]^T[\alpha\Theta_*^{-1}-\bar{H}]^{-1}[\bar{H},\bar{h}]B\right)
+\half\Tr\left(Z_2B^T[\bar{H},\bar{h}]^T[\alpha\Theta_*^{-1}+\bar{H}]^{-1}[\bar{H},\bar{h}]B\right)\leq c/\alpha
\]
for all $Z_1,Z_2\in \Z$ (recall that $\Z$ is bounded). Consequently, given $\omega>0$, we can find $\alpha=\alpha_\omega>0$ large enough to ensure that
$$
-\gamma\alpha_\omega\Theta_*^{-1}\preceq \bar{H}\preceq \gamma\alpha_\omega\Theta_*^{-1}\;\;\mbox{and}\;\; 2c/\alpha_\omega\leq\omega,
$$
which combines with (\ref{retro1}) to imply that
$$
\left[\widehat{\Psi}_+^K(\bar{h},\bar{H})+\widehat{\Psi}_-^K(\bar{h},\bar{H})\right]\leq \omega+2K^{-1}\alpha_\omega\ln(2/\epsilon),
$$
and (\ref{2onehas}) follows. \qed

\end{document}